\documentclass[12pt]{amsart}
\usepackage{amsmath,amssymb,amssymb,mathtools,hyperref}
\usepackage{microtype}
\usepackage{xcolor,listings,enumitem,needspace,fancyhdr,mathtools,booktabs,longtable}

\usepackage{geometry}
\usepackage{graphicx}
\usepackage{enumitem}
\numberwithin{equation}{section}
\theoremstyle{plain}
\newtheorem{thm}{Theorem}[section]
\newtheorem{lem}[thm]{Lemma}
\newtheorem{cor}[thm]{Corollary}

\theoremstyle{remark}
\newtheorem{rmk}[thm]{Remark}

\newcommand{\br}[1]{\left\langle#1\right\rangle}
\newcommand{\M}{\operatorname{M}}

\newcommand{\wt}{\operatorname{wt}}

\begin{document}

\title[Snowflake $q$-Flipping]{Snowflake $q$-Flipping}

\author{Tri Lai}
\address{Department of Mathematics, University of Nebraska -- Lincoln, Lincoln, NE 68505}
\email{tlai3@unl.edu}



\date{\today}

\dedicatory{}

\begin{abstract}
We prove a weighted extension of Byun’s snowflake-flipping theorem for cored hexagons with unit triangular holes along three crossing lines. The ratio of the tiling generating functions of two regions related by a snowflake flip is given by an explicit product. The result recovers Byun’s theorem when all lozenge weights are $1$ and resolves an open problem posed by the author in 2024. Our proof uses shuffling identities for semihexagons with dents on two sides.
\end{abstract}

\maketitle
\section{Introduction}\label{Intro}

Byun studied lozenge tilings of snowflake regions, obtained from a semiregular hexagon on the triangular lattice by removing a central triangle and additional unit triangles along three crossing lines. He showed that the tiling counts of two regions related by a snowflake flip differ by a simple multiplicative factor \cite[Theorem 2.1]{Byun2}. His result was inspired by earlier work of Ciucu, Rohatgi, and the author on squeezing a triad of bowties in a hexagon \cite{CLR}.  

Snowflake flipping and bowtie squeezing are instances of the tiling shuffling phenomenon introduced by Lai and Rohatgi \cite{shuffling}. Although changing the holes of a region can greatly alter its number of tilings, certain changes affect that number by a simple explicit multiplicative factor. The interesting feature is that the number-tiling ratio of the two regions related by such a rearrangement has a simple explicit expression despite the complexity of the individual tiling counts. Further examples appear in \cite{Semitwodent, shuffling2, shuffling3, Ful, FulmekShuffling2, Byun, Con, LaiByunIntrusion2, LaiByunIntrusion3}.  

In this paper, we prove a weighted extension of Byun’s result. Each vertical lozenge has weight
\[\frac{Xq^{m}+Yq^{-m}}{2}.\]
Here $X$, $Y$, and $q$ are indeterminates, and $m$ is the shifted horizontal coordinate specified in Section 2. All other lozenges have weight $1$. The weight of a tiling is the product of its lozenge weights, and the tiling generating function of a region is the sum of the weights of all its tilings. Weights of this form were inspired by the elliptic weights in \cite{Borodin} and also appear in \cite{Rosen, LaiByunIntrusion1, LaiByunIntrusion2, LaiByunIntrusion3, LR20, Semitwodent}.

Our main theorem gives an explicit product for the ratio of the tiling generating functions of two snowflake regions related by flipping. We call this phenomenon \emph{snowflake $q$-flipping}. Setting $X=Y=q=1$ makes every lozenge weight equal to $1$ and recovers Byun's theorem. The result resolves the open problem in \cite[Problem 15]{LaiOPAC}. The proof uses the author’s shuffling identities for semihexagons with dents on two sides \cite{Semitwodent}.

Section 2 defines the regions and weights and states the main theorem. Section 3 discusses connections with weighted shuffling and Rosengren’s enumeration \cite{Rosen} and the weighted shuffling theorem of Lai and Rohatgi \cite{shuffling}. It also outlines reductions for weighted extensions of the shamrock enumeration of Ciucu and Krattenthaler \cite{CK13} and the bowtie-squeezing phenomenon of Ciucu, Lai, and Rohatgi \cite{CLR}. Detailed treatments of the latter two applications will be given in \cite{LaiqShamrock} and \cite{LaiqBowtie}, respectively.

Section 4 proves the main theorem using the shuffling lemmas established in Section 5.

\section{Statement of the main theorem}\label{Statement}

Let $x,n$ be nonnegative integers, and let $H_{n,x}$ be the region obtained by removing an up-pointing triangle of side-length
$x$ from the center of the lattice hexagon with side-lengths $n+x,$ $n,$ $n+x,$ $n,$ $n+x,$ $n$  (in clockwise order starting from the northeast side\footnote{From now on, we list the side-lengths of a hexagon in this order.}). The three diagonals joining opposite vertices meet pairwise at the vertices of the central triangular hole. The six segments along these diagonals connecting a vertex of the central triangular hole to a vertex of the outer hexagon are called \textit{axes} of the region.  Label the six axes clockwise, starting with the axis joining the upper vertex of the central triangular hole to the upper-left vertex of the outer hexagon. On each axis, label the unit intervals $1,$ $2,$ \dots, $n$, starting at the central hole.

\begin{figure}[ht]\centering
 \includegraphics[width=10cm]{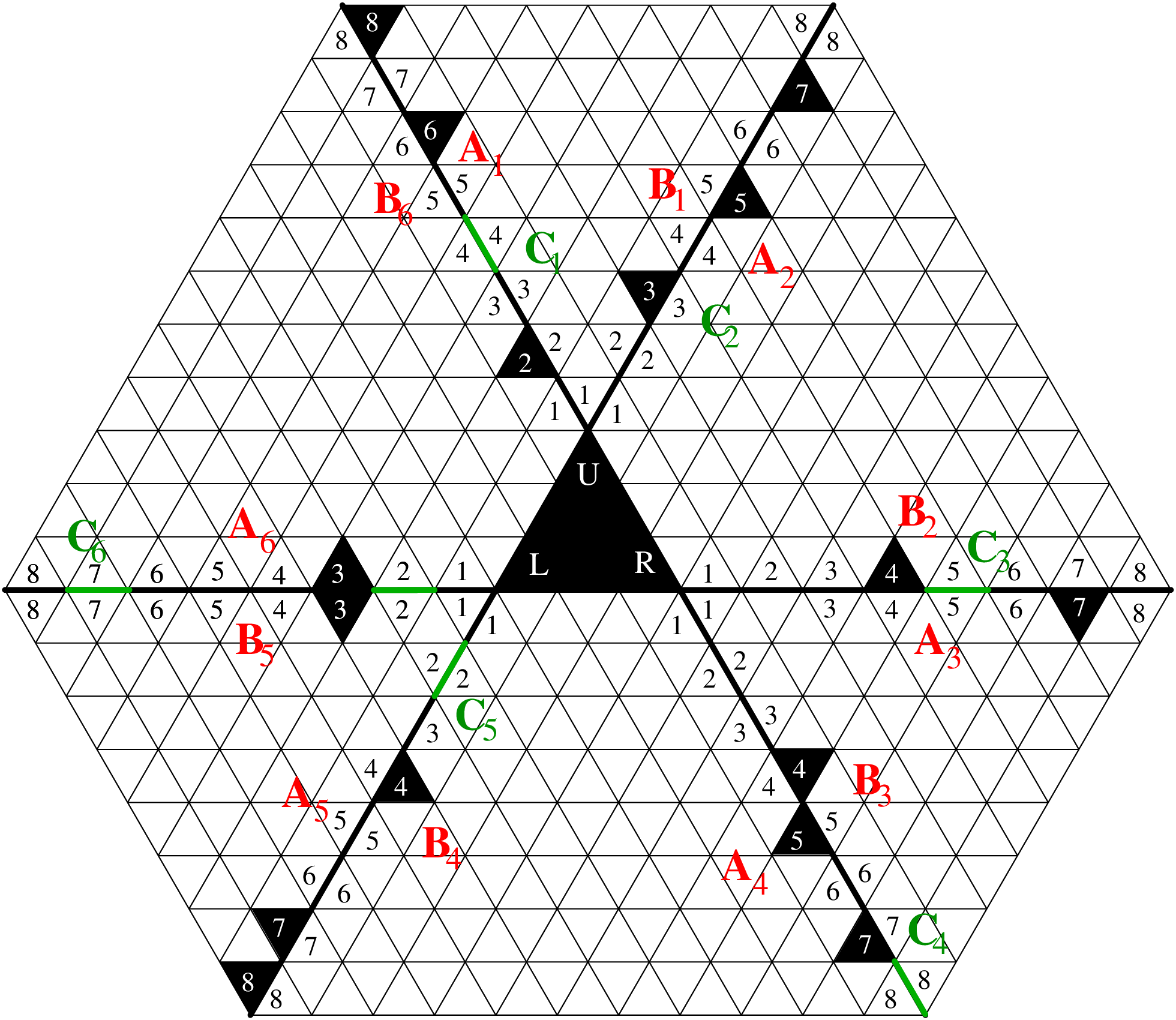}
 \caption{The snowflake region $H_{8,3}(\mathbf{A}, \mathbf{B}, \mathbf{C})$ with $A_1=\{6,8\}$, $A_2=\{5,7\}$, $A_3=\{7\}$, $A_4=\{5,7\}$, $A_5=\{7,8\}$, $A_6=\{3\}$, $B_1=\{3\}$, $B_2=\{4\}$, $B_3=\{4\}$, $B_4=\{4\}$, $B_5=\{3\}$, $B_6=\{2\}$, $C_1=\{4\}$, $C_2=\emptyset$, $C_3=\{5\}$, $C_4=\{8\}$, $C_5=\{2\}$, $C_6=\{2,7\}$. Black triangles mark the holes, and green intervals mark the barriers.}\label{snowflake}
 \end{figure}

Let $\mathbf{A}:=(A_1,A_2,\dots,A_6)$, $\mathbf{B}:=(B_1,B_2,\dots,B_6)$, and $\mathbf{C}:=(C_1,C_2,\dots,C_6)$ be $6$-tuples of  subsets of $[n]:=\{1,2,\dots,n\}$. We define $H_{n,x}(\mathbf{A}, \mathbf{B}, \mathbf{C})$ to be the region obtained from $H_{n,x}$ by removing the unit triangles touching the labelled axes at the positions specified by the sets $A_1,$ $A_2,$ $\dots,$ $A_6,$ and $B_1,$ $B_2,$ $\dots,$ $B_6$, and by placing barriers on the labelled unit intervals along the axes at the positions specified by the sets $C_1,C_2,\dots, C_6$  (see Figure \ref{snowflake} for an example).   No lozenge in an admissible tiling may cross a barrier.  In particular, going along the $i$-th axis from from the central hole toward the outer hexagon,  we remove unit triangles on the left of the axis at the position specified by the set $B_{i-1}$, remove the unit triangles on the right at the positions specified by the set $A_i$, and place barriers at the positions specified by $C_i$, for $i=1,2,\dots, 6$. Here $B_0=B_6$ by convention. In the case $n=x=0$, our region becomes the empty region with tiling generating function $1$\footnote{From now on, we always assume that an empty region has tiling generating function $1$}. In the case $n>0$ and $x=0$, the central hole collapses to the common intersection point and all six axes connect the this point to a vertex of the outer hexagon.

To make sure that the barriers do not occupy sides of the removed unit triangles along each axis, we assume that  $C_i\cap A_i=C_i\cap B_{i-1}=\emptyset$.
We also assume that $1\notin A_i\cap B_i$, for $i=1,2,\dots, 6$, to prevent a possible overlapping of  two removed triangle sets along two adjacent axes. Define  \[o:=\sum_{i=1}^{3} \left(|A_{2i-1}|+|B_{2i-1}|\right), \quad \text{and} \quad  e:=\sum_{i=1}^{3} \left(|A_{2i}|+|B_{2i}|\right),\] and assume that $o=e$. This condition is necessary for the region defined below to admit a tiling. We impose these three conditions throughout.

 
\begin{figure}[ht]\centering
  \includegraphics[width=10cm]{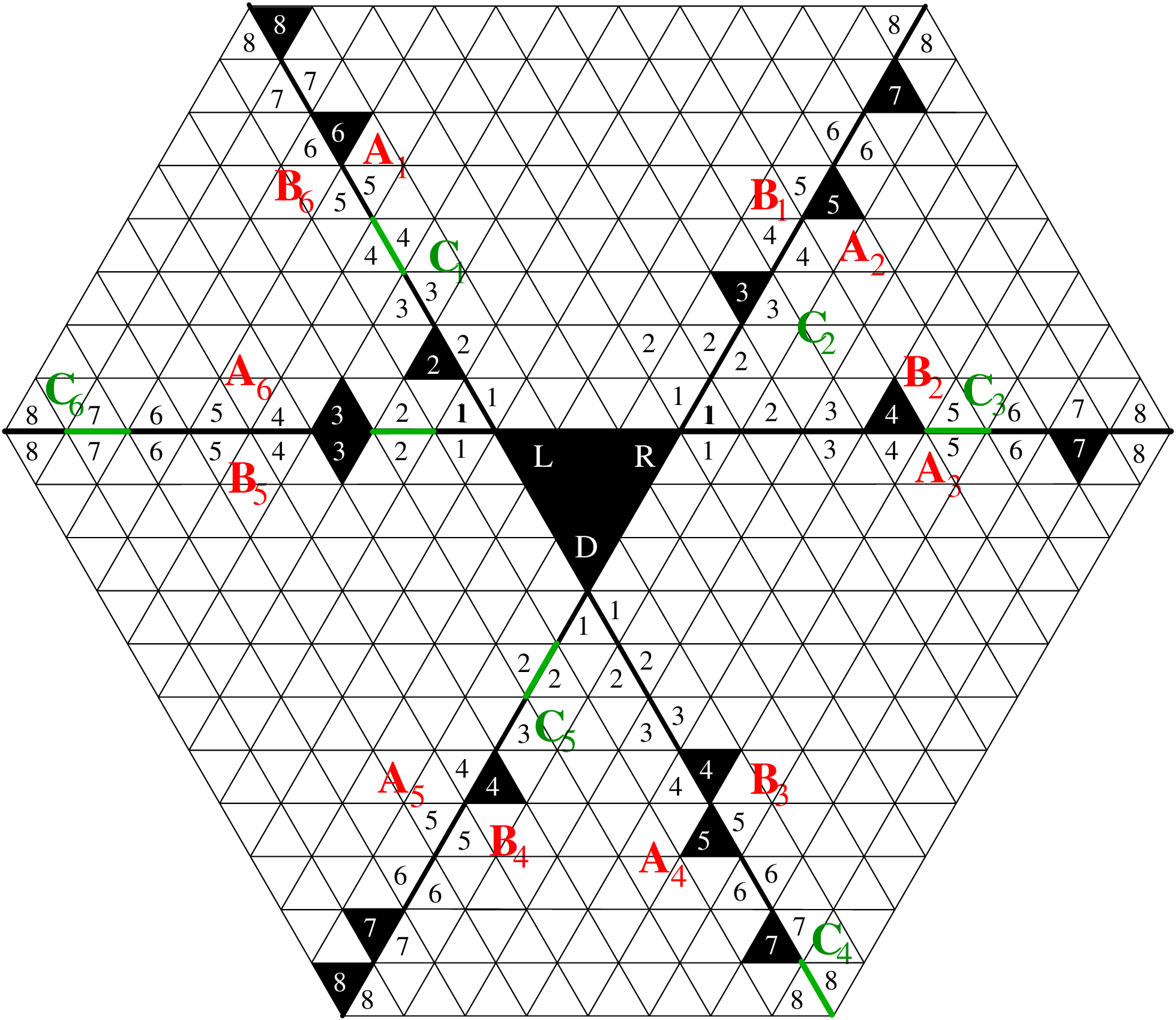}
 \caption{The snowflake region $\overline{H}_{8,3}(\mathbf{A}, \mathbf{B}, \mathbf{C})$ with the same sets $A_i$, $B_i$, and $C_i$ as in Figure \ref{snowflake}.} \label{snowflake2}
 \end{figure}

 We now define the sibling region $\overline{H}_{n,x}(\mathbf{A}, \mathbf{B}, \mathbf{C})$ of $H_{n,x}(\mathbf{A}, \mathbf{B}, \mathbf{C})$ as follows.  We start with the  hexagon $\overline{H}_{n,x}$ of side-lengths $n, n+x,n,n+x,n,n+x$ with a central down-pointing triangle of side-length $x$ removed. 
 If we also draw three main diagonals in this hexagon, then the pairwise intersections between the diagonals are exactly the vertices of the triangular hole. The six lattice segments joining vertices of the triangular hole and the vertices of the hexagon form the six axes of the new region. Starting from the axis connecting the left vertex of the triangular hole to the upper-left vertex of the outer hexagon, we list the axes by $1,2,\dots, 6$, in clockwise order.
  To obtain the region $\overline{H}_{n,x}(\mathbf{A}, \mathbf{B}, \mathbf{C})$, we remove the unit triangles at the positions specified by $A_1,A_2,\dots,A_6,$ and $B_1,B_2,\dots, B_6$, and place the barriers at the positions specified by $C_1,C_2,\dots, C_6$. We remove unit triangles and place the barriers at the same labelled positions along corresponding axes in the two regions.  See Figure \ref{snowflake2} for an example.
 
\begin{figure}[ht]\centering
     \includegraphics[width=7cm]{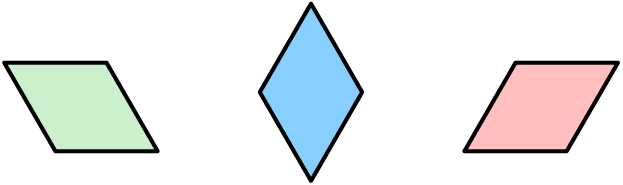}
 \caption{From left to right, the three lozenge orientations are left-tilting,
vertical, and right-tilting.} \label{lozengeorient}
 \end{figure}

\begin{figure}[ht]\centering
   \includegraphics[width=15cm]{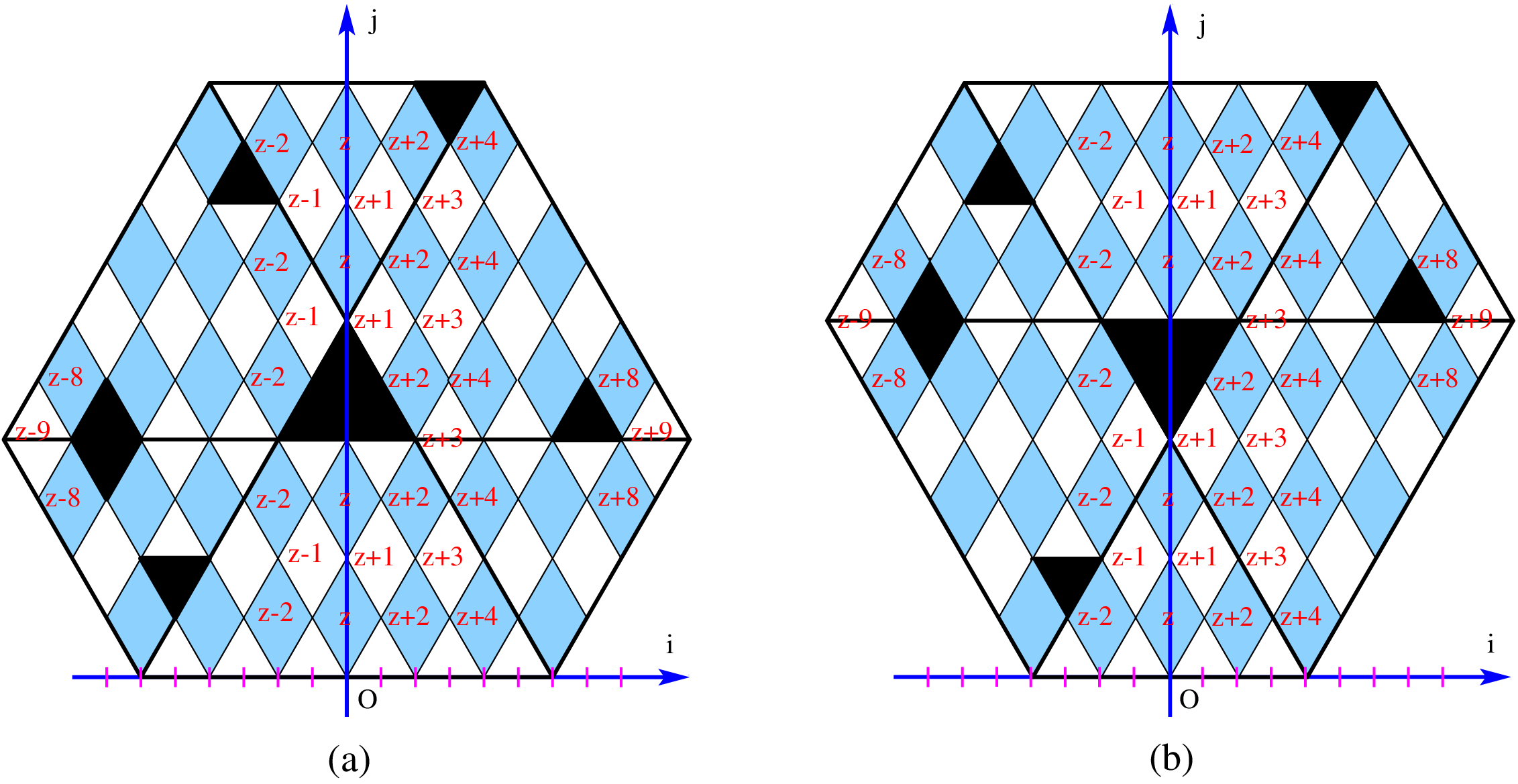}
 \caption{The weight assignment $\operatorname{wt}_{X,Y,z}$ for (a)  the snowflake region $H_{n,x}(\mathbf{A},\mathbf{B},\mathbf{C})$ and (b) its sibling region $\overline{H}_{n,x}(\mathbf{A},\mathbf{B},\mathbf{C})$.  Each labelled vertical lozenge has weight
$f_{X,Y}(k)=\frac{Xq^{k}+Yq^{-k}}{2}$ when its displayed label is $k$; all other lozenges have weight $1$.} \label{snowweight}
 \end{figure}

We now define the weight assignment $\operatorname{wt}_{X,Y,z}=\operatorname{wt}_{X,Y,z}(q)$ in our snowflake regions, where $X,Y,q$ are indeterminates and $z$ is an integer. We introduce an orthogonal coordinate system $Oij$ on the  triangular lattice so that the horizontal $i$-axis runs along the base of the region $H_{n,x}(\mathbf{A}, \mathbf{B}, \mathbf{C})$, and the vertical $j$-axis passes through the top vertex of the central triangular hole. One unit on the $i$-axis is half the side length of a unit triangle, and one unit on the $j$-axis is its altitude. The weights depend only on the $i$-coordinate (the choice of unit on the $j$-axis does not matter as our weighting does not depend on the $j$-coordinate).  A vertical
lozenge whose center has $i$-coordinate $k$, for some integer $k$, has weight
 \[f_{X,Y}(k+z)=\frac{Xq^{k+z}+Yq^{-k-z}}{2}.\]
Figure \ref{lozengeorient} shows the three orientations of lozenges. In particular, every vertical
lozenge centered on the $j$-axis has weight $f_{X,Y}(z)$. All lozenges of other two orientations have weight $1$. See Figure \ref{snowweight}(a) for an example of the weight assignment for lozenges of the region $H_{n,x}(\mathbf{A}, \mathbf{B}, \mathbf{C})$. The \emph{weight} of a tiling $\tau$  is the product of the weights of its lozenges. For a region $R$, we define
\[\M_{X,Y,z}\left( R\right) =\sum_{\tau}\operatorname{wt}_{X,Y,z}(\tau),\]
where the sum is taken over all tilings $\tau$ of $R$\footnote{We always assume in this paper that a region $R$ admits a tiling (or is tileable) if and only if its tiling generating function $\M_{X,Y,z}\left( R\right)$ is different from $0$.}. We call $\M_{X,Y,z}\left( R\right)$ the \textit{tiling-generating function} of $R$ and suppress the parameter $q$ from the notation when no confusion can arise. We have defined the generating function $\M_{X,Y,z}\left( H_{n,x}(\mathbf{A}, \mathbf{B}, \mathbf{C})\right)$.

For the sibling region $\overline{H}_{n,x}(\mathbf{A},\mathbf{B},\mathbf{C})$, use the same rule, with the vertical $j$-axis passing through the bottom vertex of its central triangular hole (see Figure \ref{snowweight}(b)).


For fixed $X,$ $Y,$ $z,$ and $q,$ the weights depend on the chosen horizontal coordinates. 

 \bigskip
 
For any positive integer $n$, define the \emph{symmetric $q$-number} by \[\br{n}=\br{n}_q:=\frac{q^{n}-q^{-n}}{q-q^{-1}},\]  and set $\br{0}=0$. As $\lim_{q\to1}\br{n}_q=n$, we define its value at $q=1$ to be $n$.
For a nonnegative integer $x$, define the  \emph{symmetric $q$-Pochhammer symbol} by \[\br{n}_x=\br{n}_{x,q}:=\prod_{i=0}^{x-1} \br{n+i},\] with the empty product equal to $1$.

We define two operators  \[E_{x,z}(S)=\prod_{i=1}^{x} \prod_{s\in S} f_{X,Y}(z+s-i),\] 
 \[E'_{x,z}(S)=\prod_{i=1}^{x} \prod_{s\in S} f_{Y,X}(-z+s-i),\] where  $S$ is a finite set of integers and empty products are equal to $1$ by convention. Denote by
 $A+k:=\{a+k \mid a\in A\}$  the  translation of $A$ by $k$.

 
\begin{thm}\label{mainthm} Let $n$ and $x$ be nonnegative integers, let $z$ be an integer, and let $X,$ $Y,$ and $q$ be indeterminates. Suppose  $A_1, $ $A_2,$ $\dots,$ $A_6,$ $B_1,$ $B_2,$ $\dots,$ $B_6,$ $C_1,$ $C_2,$ $\cdots,$  $C_6$ are subsets of $[n]:=\{1,2,\dots,n\}$ such that $C_i\cap A_i=C_i\cap B_{i-1}=\emptyset$ and $1\notin A_i \cap B_i$,  for $i=1,2,\dots, 6$ and $B_0=B_6$ by convention.  If $\overline{H}_{n,x}(\mathbf{A},\mathbf{B},\mathbf{C})$ has a tiling, then
\begin{align}\label{maineq}
\frac{\M_{X,Y,z}(\overline{H}_{n,x}(\mathbf{A},\mathbf{B},\mathbf{C}))}{\M_{X,Y,z}(H_{n,x}(\mathbf{A},\mathbf{B},\mathbf{C}))}&=\frac{E_{x,z}(A_3+x)E'_{x,z}(B_5+x)}{E_{x,z}(B_2+x)E'_{x,z}(A_6+x)} \prod_{j=1}^{3}\frac{\prod_{i\in A_{2j-1}}\br{i}_x \prod_{i\in B_{2j-1}}\br{i}_x}{\prod_{i\in A_{2j}}\br{i}_x \prod_{i\in B_{2j}}\br{i}_x }.
\end{align}
\end{thm}

By letting $C_1=\cdots= C_6=\emptyset$ and $X=Y=q=1$, the generating functions count exactly the tilings of the two regions in the theorem. In this specialization, our symmetric $q$-integer $\br{n}$ becomes the regular integer $n$, and $f_{X,Y}(k)=1$ for any integer $k$, so $E_{x,z}(A_3+x)=E'_{x,z}(B_5+x)=E_{x,z}(B_2+x)=E'_{x,z}(A_6+x)=1$. Then our theorem implies Byun's  \cite[Theorem 2.1]{Byun2}.

The expression on the right-hand side of (\ref{maineq}) is independent of the barrier sets $C_1,\cdots, C_6$. Thus these sets can affect the two generating functions individually while leaving their ratio unchanged, whenever the regions satisfy the theorem’s hypotheses. We will showed later in the proof the theorem in Section 4 that the region $H_{n,x}(\mathbf{A},\mathbf{B},\mathbf{C})$ admits a tiling if and only if the region $\overline{H}_{n,x}(\mathbf{A},\mathbf{B},\mathbf{C})$ admits a tiling.

\begin{rmk} 1) 
Here is a short explanation for the ``balancing" condition
\[o=\sum_{i=1}^{3} |A_{2i-1}|+|B_{2i-1}|=\sum_{i=1}^{3} |A_{2i}|+|B_{2i}|=e\]
 A necessary condition for a region to admit a lozenge tiling is that it contain equally many
up-pointing and down-pointing unit triangles, since each lozenge covers one triangle of each
orientation. After the central triangular hole is removed from the initial hexagon in the definition of $H_{n,x}(\mathbf{A},\mathbf{B},\mathbf{C})$, the region is balanced. The additional removals delete $o$ down-pointing unit triangles and $e$ up-pointing unit triangles, so balance requires $o=e$.  

2) We also note that, in our setup, we allow barriers at the positions specified by the sets $C_1,C_2,\dots,C_6$. Barriers also
appear in the shuffling theorem of Lai and Rohatgi \cite[Theorem 2.4]{shuffling}.

\end{rmk}

 \begin{figure}[ht]\centering
 \includegraphics[width=15cm]{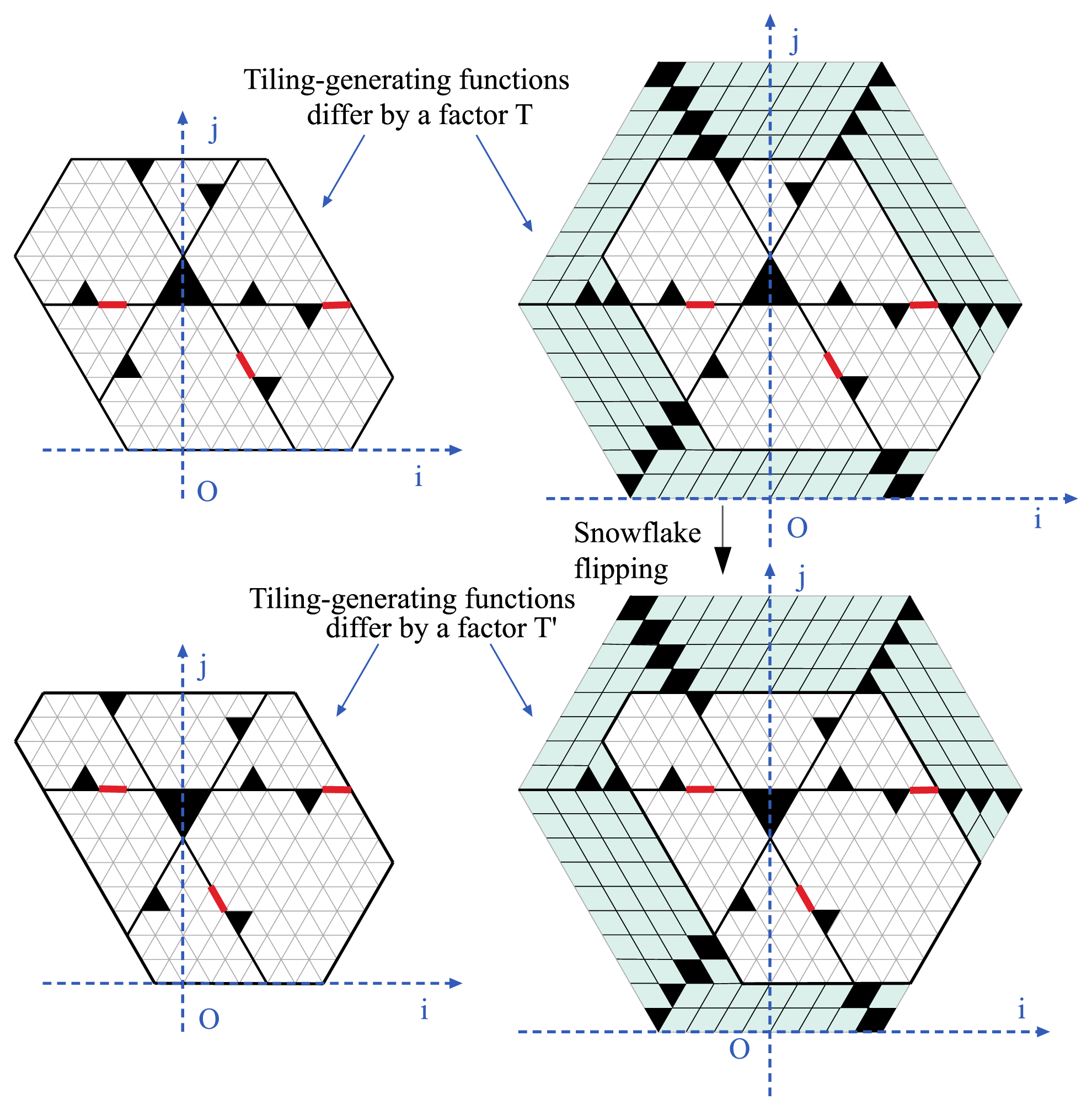}
 \caption{Reduction of snowflake $q$-flipping for a general hexagon with holes along three
crossing lines to Theorem \ref{mainthm}. Black triangles indicate holes, red intervals indicate barriers,
and shaded lozenges are forced. The additional holes are chosen so that the two products of
forced-lozenge weights, $T$ and $T'$, are equal.}\label{GneralFlipPic}
 \end{figure}

\section{Applications to weighted enumeration of tilings}

This section outlines applications and reductions arising from Theorem \ref{mainthm}; detailed treatments of the shamrock and bowtie weighted enumerations will appear in \cite{LaiqShamrock} and \cite{LaiqBowtie}.

\subsection{A reduction for general hexagons and Rosengren Corollary 2.2 \cite{Rosen}}

The extension of snowflake $q$-flipping to more general hexagons can be approached using the
construction in \cite{Twofern}, also discussed in \cite{Byun2}. A complete treatment requires a careful analysis of
the possible configurations. 
Figure  \ref{GneralFlipPic} illustrates one of the confgurations.

In each row of Figure \ref{GneralFlipPic}, the tiling generating function of the region on the right equals that of
the region on the left multiplied by the product of the weights of the shaded forced lozenges. Let the left-hand regions in the two rows be $L$ and $L'$, the right-hand snowflake regions $S$ and $S'$, and the forced-lozenge products $T$ and $T'$. Then the stated reduction gives
\[\frac{M_{X,Y,z}(L)}{M_{X,Y,z}(L')}=\frac{T'}{T} \frac{M_{X,Y,z}(S)}{M_{X,Y,z}(S')}.\]
Choose the additional holes so that the two multiplicative factors $T$ and $T'$ are equal. The ratio of the
generating functions for the two regions $L$ and $L'$ on the left then equals the ratio for the two snowflake
regions $S$ and $S'$ on the right, to which Theorem \ref{mainthm} applies.

Now if we consider the case with no barrier and no unit triangle removed along the axes, then our snowflake regions are simply cored hexagons. The snowflake $q$-flipping only flips the central triangular holes, and our theorem (with $X=Y=1$ and $z=0$) explains why the tiling generating functions of the two cored hexagons in Rosengren Corollary 2.2 \cite{Rosen} differ by a simple product.

\subsection{A weighted shuffling phenomenon}

\begin{figure}[ht]\centering
  \includegraphics[width=12cm]{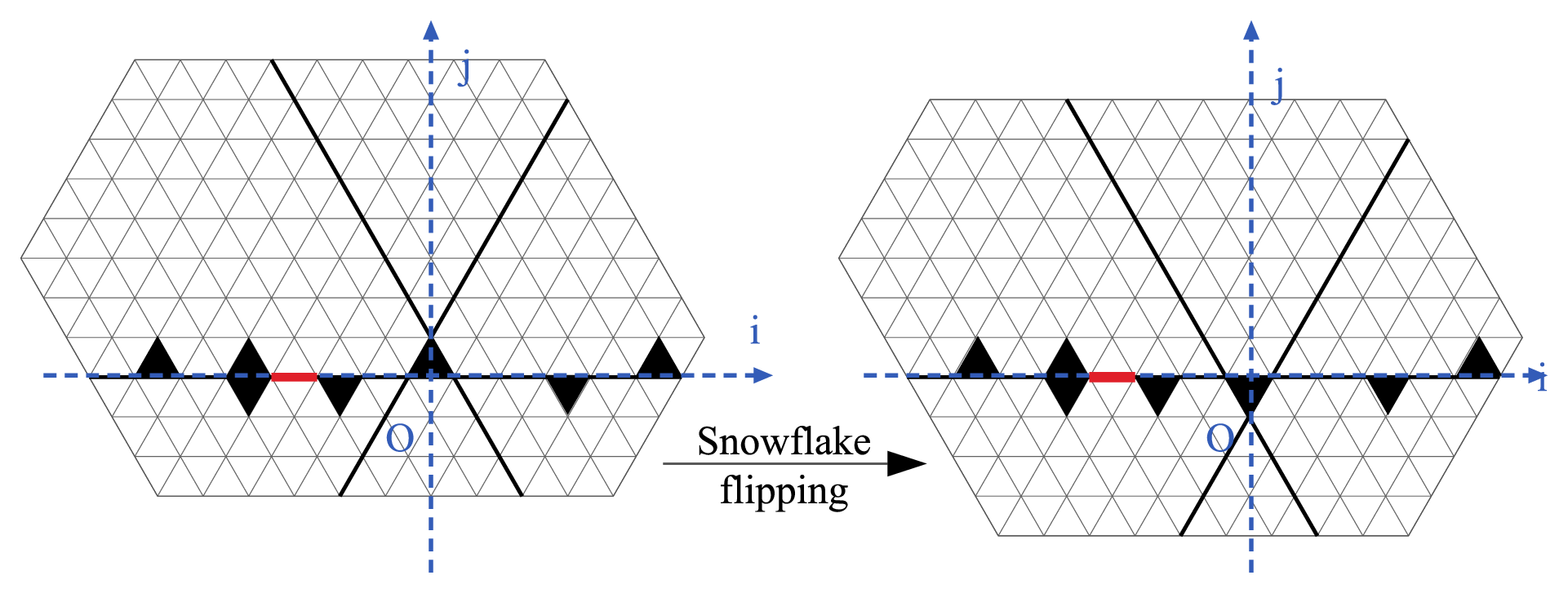}
 \caption{A snowflake $q$-flip performs
one shuffling step.}\label{ShufflingPic}
 \end{figure}

Applying Theorem \ref{mainthm} to the pair of regions in Figure \ref{ShufflingPic} gives one shuffling step. Here we assume that there are no barriers, except for the ones along the horizontal lines. Repeating this operation yields conceptually the $q$-shuffling theorem of Lai and Rohatgi  \cite[Theorem 2.4]{shuffling}.
Byun’s unweighted snowflake-flipping theorem gives the unweighted specialization of the
Lai--Rohatgi theorem, namely, Theorem 2.3 in \cite{shuffling}.  Further
applications of the generalized shuffling phenomenon are planned for  \cite{Lainewshuffling}.

 \begin{figure}[ht]\centering
   \includegraphics[width=15cm]{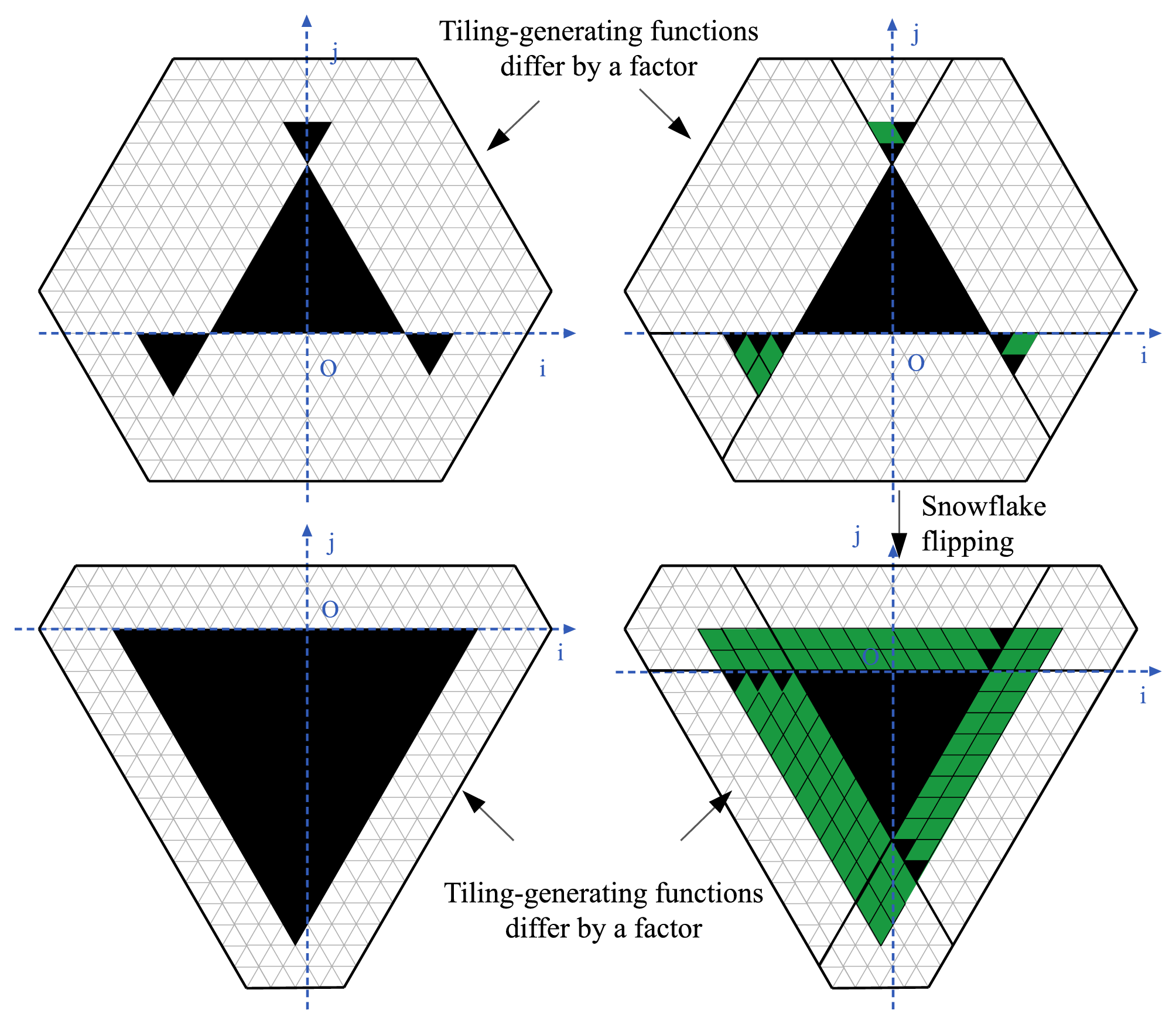}
 \caption{Reduction of the weighted enumeration of a hexagon with a shamrock hole (upper-left) to that of
a cored hexagon (lower-left), connecting the problem with Rosengren’s enumeration in \cite{Rosen}.  Black triangles indicate holes
and shaded lozenges are forced. In each row, the generating functions of the two regions differ by a multiplicative factor equal to the product of the weights of the forced lozenges. The two regions in the second column are related by a snowflake $q$-flip.}\label{ShamrockPic}
 \end{figure}

\subsection{Hexagons with a shamrock hole}

Theorem \ref{mainthm} (with $C_i=\emptyset$ for any $i$) reduces the weighted
enumeration of lozenge tilings of a hexagon with a central shamrock hole to the corresponding problem
for a cored hexagon as illustrated in Figure \ref{ShamrockPic}.  Theorem \ref{mainthm} (with $X=Y=1$ and $z=0$) expresses the ratio of these generating functions as an explicit product.  Using Rosengren’s enumeration of the region with a triangular hole \cite{Rosen}, we obtain a product formula for the tiling generating function of a hexagon with a shamrock hole. This extends the result of  Ciucu and Krattenthaler \cite{CK13}. 

\subsection{A $q$-analogue of bowtie squeezing}

Theorem \ref{mainthm} (with $C_i=\emptyset$ for any $i$) also provides a route to a $q$-analogue of the bowtie-squeezing phenomenon of
Ciucu, Lai, and Rohatgi \cite[Theorem 1]{CLR}. We seek the ratio of the tiling generating functions of the two regions $R$ and $Q$ that are related by a sequence of bowtie squeezes. To do so, we compare their tiling generating functions with that of an intermediate region $\overline{R}$. Here $\overline{R}$ is obtained from $R$ and also from $Q$ by completely squeezing all three bowties. (The sums of sizes of the two lobes in the corresponding bowties of $R$ and $Q$ are the same, so their completely-squeezed version are the same, which is the region $\overline{R}$.)

\begin{figure}[ht]\centering
   \includegraphics[width=15cm]{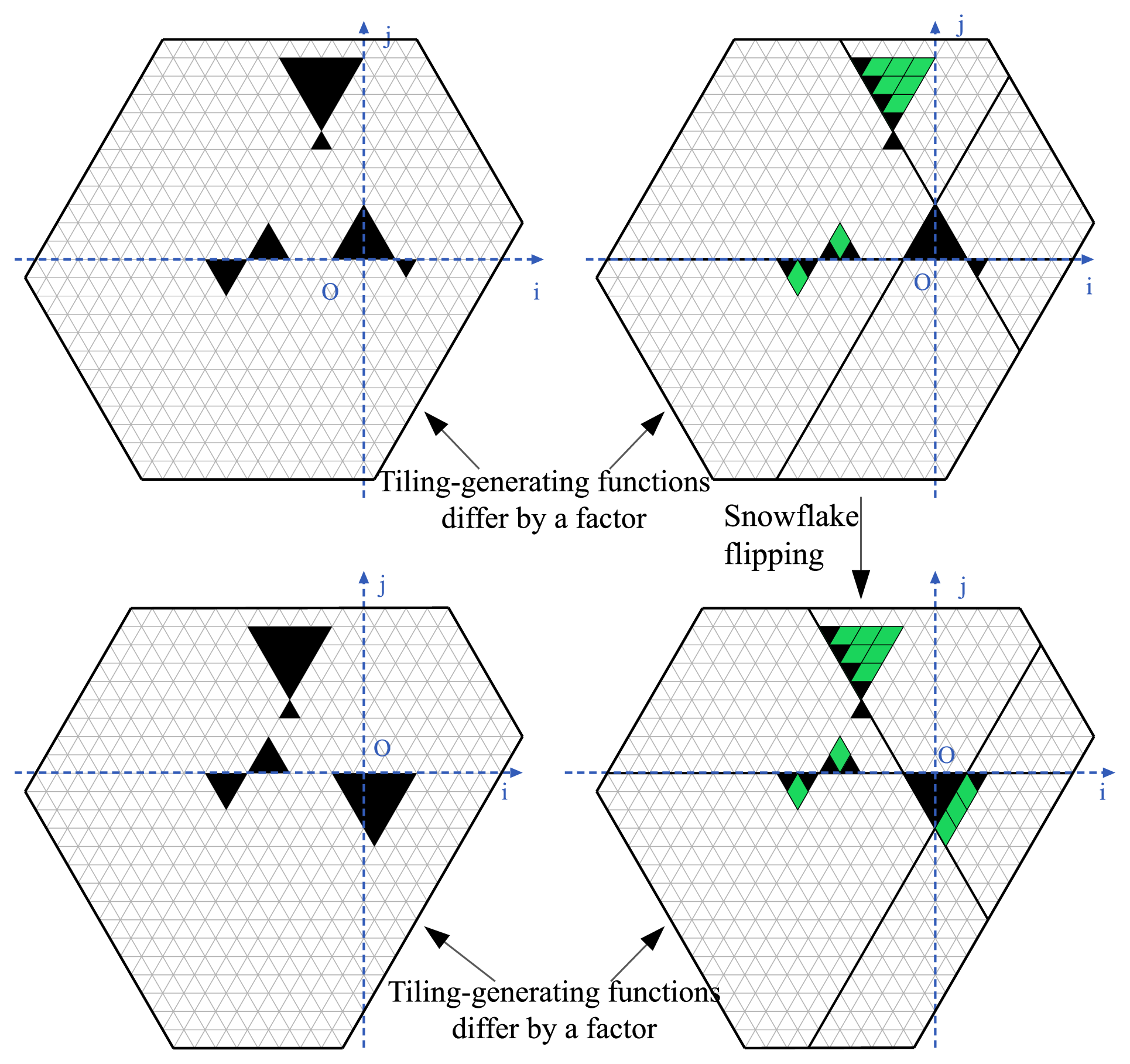}
 \caption{Obtaining  a weighted enumeration version of bowtie squeezing by Ciucu, Lai, and Rohatgi \cite{CLR}. Completely squeezing one bowtie in the upper-left region yields the lower-left region. Black triangles indicate holes
and shaded lozenges are forced. In each row, the generating functions of the two regions differ by a multiplicative factor equal to the product of the weights of the forced lozenges. The two regions in the second column are related by a snowflake $q$-flip.}\label{BowtiePic}
 \end{figure}

 Each application squeezes one bowtie of $R$ completely; see Figure \ref{BowtiePic}.
Applying the theorem to the other two bowties gives a product expression for the ratio of the
tiling generating functions of $R$ and $\overline{R}$. Repeating the argument for $Q$
and $\overline{R}$ gives a second ratio. Dividing $\frac{\M_{X,Y,z}(R)}{\M_{X,Y,z}(\overline{R})}$ by $\frac{\M_{X,Y,z}(Q)}{\M_{X,Y,z}(\overline{R})}$ yields $\frac{\M_{X,Y,z}(R)}{\M_{X,Y,z}(Q)}$. 

\section{Proof of the main theorem}

The proof proceeds by decomposing each snowflake region into three component regions after fixing its crossing lozenges. We first establish the weight transformations and semihexagon identities needed to compare the corresponding components. We then show that the resulting ratio is independent of the crossing positions and sum over all compatible choices.

We first investigate some useful properties of our weight assignments.

\begin{enumerate}
 
\item[(P1)] The weight function $f_{X,Y}$ satisfies the identity:
 \[f_{X,Y}(n)=\frac{Xq^{n}+Yq^{-n}}{2}= \frac{Y \cdot q^{(-n)}+Xq^{-(-n)}}{2}=f_{Y,X}(-n),\]
 for any integer $n$. With the coordinate system held fixed, the identity gives a useful reflection
symmetry. Let $R^{|}$ be the reflection of $R$ over the vertical $j$-axis, and let $\phi$ map each tiling of $R^{|}$ to its reflected tiling of $R$. A vertical lozenge at horizontal coordinate $k$ in $R^{|}$ has weight $f_{X,Y}(k+z)$. It image has horizontal coordinate $-k$. With parameters $(Y,X,-z)$, its weight is
 \[f_{Y,X}(-k-z)=f_{X,Y}(k+z).\]
 Thus $\wt_{X,Y,z}(\tau)=\wt_{Y,X,-z}(\phi(\tau)).$
Summing over all tilings $\tau$ of $R^{|}$ yields
 \begin{equation}\label{vreflectprinciple}
 \M_{X,Y,z}(R^{|})= \M_{Y,X,-z}(R).
 \end{equation}

\item[(P2)]  With the coordinate system held fixed, reflection over the horizontal $i$-axis preserves the
horizontal coordinate of every lozenge and therefore preserves its weight. If $R^{-}$ is the reflection of $R$ over this axis, then
 \begin{equation}\label{hreflectprinciple}
 \M_{X,Y,z}(R^{-})= \M_{X,Y,z}(R).
 \end{equation}
 Similarly, if $R'$ is obtained from $R$ by a vertical translation, then \[\M_{X,Y,z}(R')=\M_{X,Y,z}(R).\]
 
\item[(P3)]  With the coordinate system held fixed, let $R^{\to}_{l}$ denote the translate of $R$ by $l$ units along the horizontal  axis, where $l$ is an integer. Here if $l$ is positive  we translate to the right, and if $l$ is negative, we translate to the  left by $-l$ units.
 Let $\psi$ map each tiling of $R^{\to}_{l}$ back to  a tiling of $R$ by by translating it $-l$ units in the horizontal direction.  A vertical lozenge at horizontal coordinate $k$ in $R^{|}$ has weight $f_{X,Y}(k+z)$. Its image under $\psi$ has horizontal coordinate $k-l$ and, with the parameters $(X,Y,z+l)$, has weight \[f_{X,Y}(k-l +(z+l))=f_{X,Y}(k+z).\]  The weights of the other two lozenge orientations remain $1$. Thus $\wt_{X,Y,z}(\tau)=\wt_{X,Y, z+l}(\psi(\tau)).$
 Summing over all tilings of $R^{\to}_{l}$ yields
 \begin{equation}\label{vreflectprinciple}
 \M_{X,Y,z}(R^{\to}_{l})= \M_{X,Y,z+l}(R).
 \end{equation}
\end{enumerate}

\bigskip

We next introduce the semihexagonal regions used in the proof of the main theorem.

Let $n,x$ be two nonnegative integers. Let $S_{n,x}$ be the upper half of a  hexagon with side-lengths $n,x,n,n,x,n$. The region is a trapezoid (when $x=0$, the trapezoid degenerates to a triangle). We label the unit triangles along its two oblique sides from top to bottom by $1,2,\dots, n$. For subsets $A$ and $B$  of $[n]$,  remove the up-pointing unit triangles at the positions specified by $A$ on the left side and $B$ on the right side. We denote the resulting region by $S_{n,x}(A,B)$. We assume that $|A|+|B|=n$, since otherwise the region does not admit any tiling. 

\begin{figure}[ht]\centering
  \includegraphics[width=15cm]{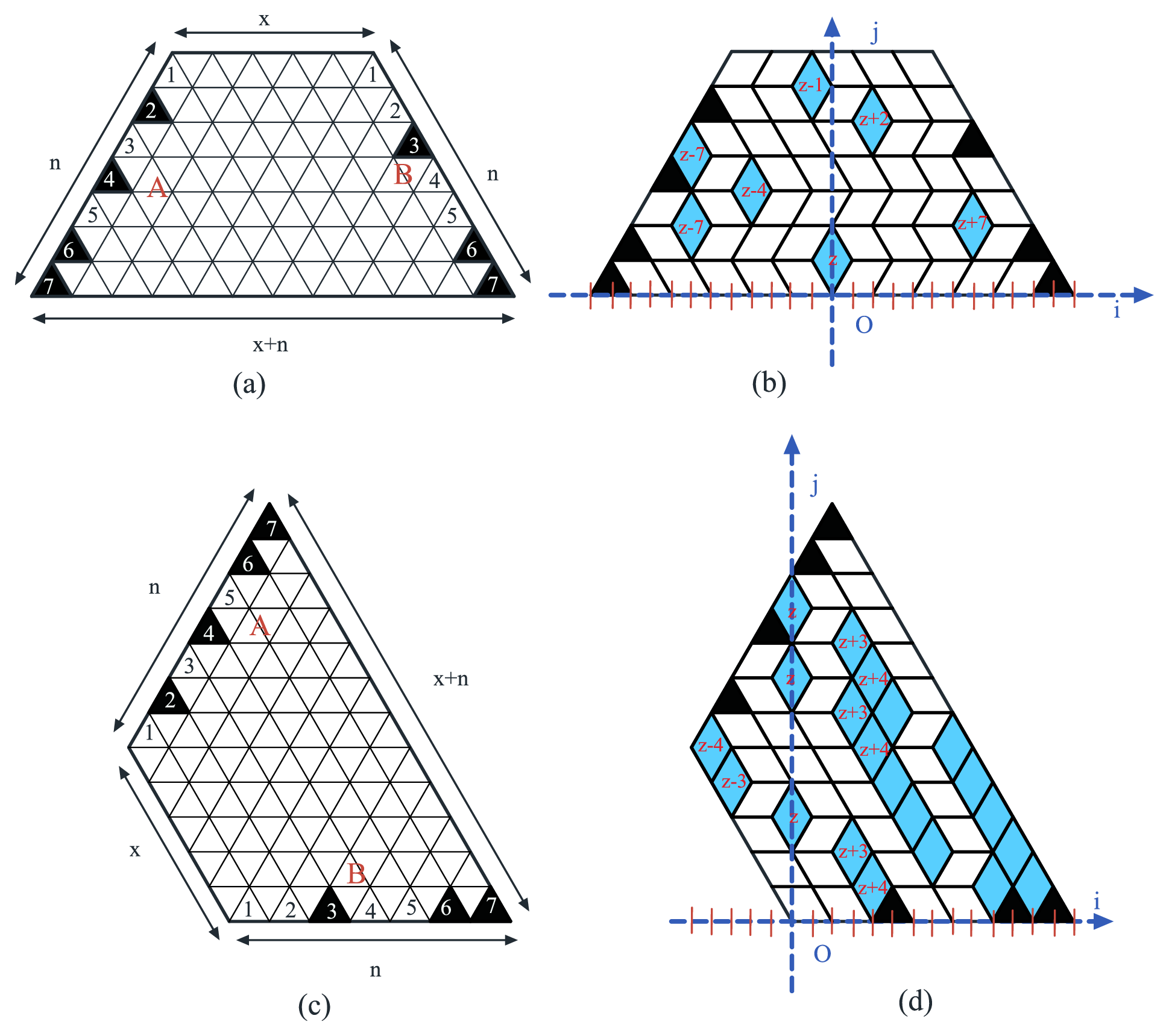}
 \caption{The semihexagons with unit triangles removed from both sides $S_{7,5}(A,B)$ in (a) and $S'_{7,5}(A,B)$ in (c) with $A=\{2,4,6,7\}$ and $B=\{3,6,7\}$ and sample tilings in (b) and (d), respectively. Each labelled vertical lozenge has weight
$f_{X,Y}(k)=\frac{Xq^{k}+Yq^{-k}}{2}$ when its displayed label is $k$; all other lozenges have weight $1$.}\label{Sshape}
 \end{figure}

To assign weights, we place an orthogonal coordinate system $Oij$ so that the horizontal $i$-axis runs along the base of the region and that the vertical $j$-axis passes through the midpoints of the two parallel sides of the region. Apply the weight assignment $\wt_{X,Y,z}$. Figure \ref{Sshape}(a) shows the region and  Figure \ref{Sshape}(b) shows a sample tiling.

Let $S'_{n,x}(A,B)$ be the southwest half of a hexagon with side-lengths $n,n,x,n,n,x$. The region is again a trapezoid (rotated by $120^{\circ}$). We label the unit triangles along its base and its northwest side by $1,2,\dots, n$ as shown in (see Figure  \ref{Sshape}(c)). For subsets $A$ and $B$  of $[n]$,  remove the up-pointing unit triangles at the positions specified by $A$ on the northwest side and $B$ on the based. 
To determine the weight assignment of the resulting region,  we place a new coordinate system so that the $i$-axis runs along the base of the resulting region and that the $j$-axis passes through the left endpoint of that base (see Figure  \ref{Sshape}(d)).

Under the standing cardinality assumption $|A|+|B|=n$, the following lemma gives a necessary and sufficient condition for the regions  $S_{n,x}(A,B)$ and  $S'_{n,x}(A,B)$ to admit a tiling.
 
 \begin{lem}[Lemma 2.1 in \cite{Semitwodent}]\label{lem0} Assume that $n,x$ are nonnegative integers, and that $A$ and $B$ are two subsets of $[n]$ with $|A|+|B|=n$. The two regions
 $S_{n,x}(A,B)$ and  $S'_{n,x}(A,B)$ have a tiling if and only if
 \[ |A\cap[j] |+|B\cap [j]|\leq t,\]
 for every integer $j$ with $1\le j\leq n$. \end{lem}
A pair $(A,B)$ satisfying the hypotheses and inequalities of Lemma \ref{lem0} is called \emph{compatible}.


For nonnegative integers $x,r$ and a finite set of integers $T$, define
 \[D_x(T)=\prod_{t\in T}\br{t}_x \quad \text{and} \quad F_x(r)=D_{x}([r])=\prod_{i=1}^{r}\br{i}_x,\]
 where empty products are equal to $1$.
In particular, $D_0(T)=F_0(r)=1$.
 We will use the following simplification for the our next proofs.
We have the following shuffling phenomenon for our regions $S_{n,x}(A,B)$ and $S'_{n,x}(A,B)$.

\begin{lem}\label{shufflingthm1} Let
$n,$ $x,$ $y$ be nonnegative integers. Let $(A, B)$ be a compatible pair of subsets of $[n]$. If $S_{n,y}(A,B)$ has a tiling, then
\begin{align*}
\frac{\M_{X,Y,z}(S_{n,x}(A,B))}{\M_{X,Y,z}(S_{n,y}(A,B))}=\frac{F_y(n)}{F_x(n)} \frac{D_x(A)D_x(B)}{D_y(A)D_y(B)}.
\end{align*}
\end{lem}

\begin{lem}\label{shufflingthm2} Under the assumptions of Lemma \ref{shufflingthm1}, if $S'_{n,y}(A,B)$ has a tiling, then
\begin{align*}
\frac{\M_{X,Y,z}(S'_{n,x}(A,B))}{\M_{X,Y,z}(S'_{n,y}(A,B))}=& \frac{E_{x,z}(B)}{E_{y,z}(B)} \cdot \frac{F_y(n)}{F_x(n)} \frac{D_x(A)D_x(B)}{D_y(A)D_y(B)}.
\end{align*}
\end{lem}

Lemmas \ref{shufflingthm1} and \ref{shufflingthm2} generalize the two shuffling results in \cite[Theorems 1.1 and 1.2]{Semitwodent}, respectively. We express our formulas in terms of the symmetric $q$-number ($\br{n}=\frac{q^n-q^{-n}}{q-q^{-1}}$), whereas \cite{Semitwodent} uses the ordinary $q$-number ($[n]_q=1+q+\cdots+q^{n-1}$).

Lemma \ref{shufflingthm1} follows from \cite[Theorems 1.1]{Semitwodent} by the translation principle (P3). The proof of that theorem also applies directly, so we omit the details. Fulmek \cite{FulmekShuffling2} gives an alternative proof of \cite[Theorems 1.1]{Semitwodent}. The formula in the original statement of \cite[Theorems 1.2]{Semitwodent} contained a typographical error. For clarity, we give a detailed proof of Lemma \ref{shufflingthm2} in Section 5.



We will use Lemmas \ref{shufflingthm1} and \ref{shufflingthm2} to prove the following weighted versions of Lemma 3.1 in \cite{Byun2}.

  \begin{figure}[ht]\centering
   \includegraphics[width=15cm]{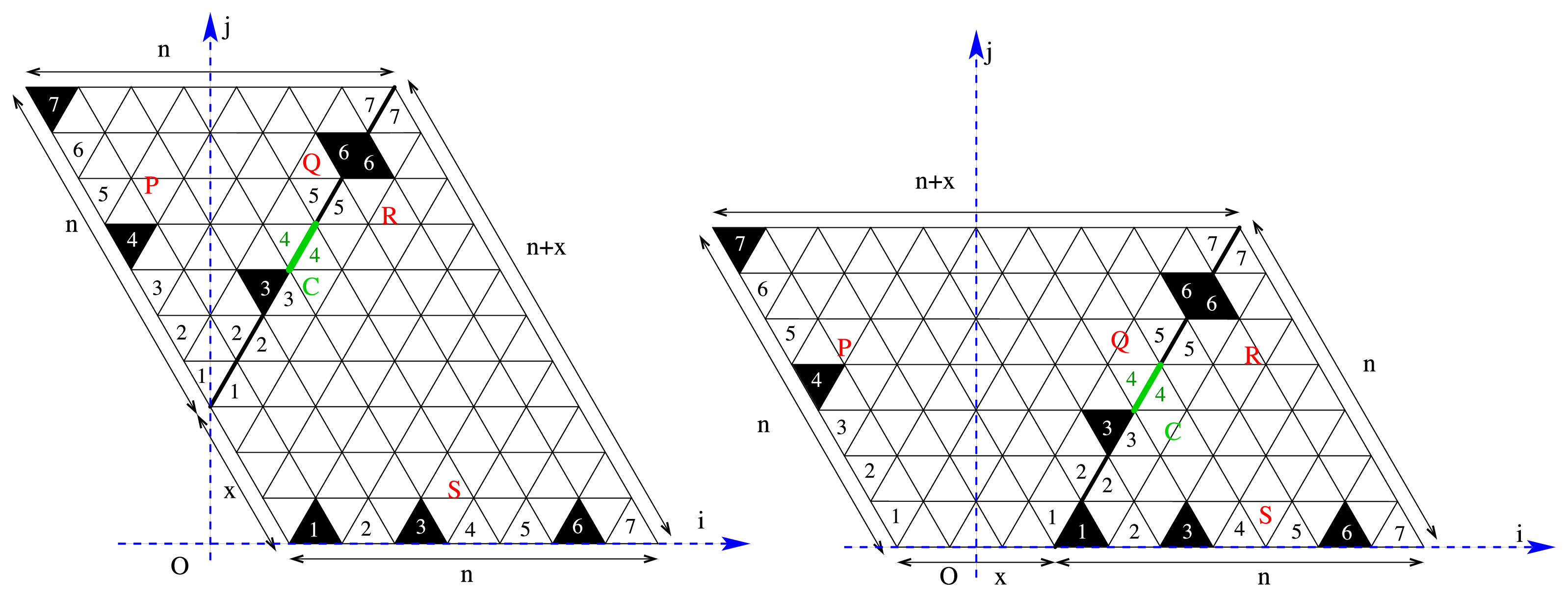}
 \caption{The regions $L^{(1)}_{7,3}(P,Q,R,S,C)$ (left) and  $\overline{L}^{(1)}_{7,3}(P,Q,R,S,C)$ (right) with $P=\{4,7\}$, $Q=\{3,6\}$, $R=\{6\}$, $S=\{1,3,6\}$, and $C=\{4\}$. Black triangles indicate removed unit triangles, and green segments
indicate barriers. The dashed blue lines show the coordinate axes.} \label{Lshape1}
 \end{figure}

Let $L^{(1)}_{n,x}$ be a leftward-slanting parallelogram of side-lengths $n,n+x,n,n+x$ (in clockwise order, starting from the top side\footnote{From now on, unless indicated otherwise, we always list the side-lengths of a parallelogram in this order.}). We divide the parallelogram into two components by a diagonal cut from the upper-right vertex to the left side. The upper component is a down-pointing triangle of side-length $n$; the lower component is a trapezoid. We label the unit positions along the indicated sides of the triangle and trapezoid by $1,2,\dots,n$, as in the left picture in Figure \ref{Lshape1}. Let $P,Q,R,S, C$ be five subsets of $[n]$.  Remove the down-pointing unit triangles indexed by $P$ and $Q$ from the triangular component and the
up-pointing unit triangles indexed by $R$ and $S$ from the trapezoidal component.  Finally, we place barriers at the positions on the cut indexed by $C$; no lozenge may cross a barrier (see the left picture in  Figure \ref{Lshape1} for an example). We denote by $L^{(1)}=L^{(1)}_{n,x}(P, Q, R, S,C)$ the resulting region. To keep the barriers disjoint from the holes along the cut, we assume that $C\cap (Q\cup R)=\emptyset$. We also assume that  $|P|+|Q|=|R|+|S|\leq n$, $1\notin P\cap Q$, and $1\notin R\cap S$ (otherwise the region does not admit a tiling).

To specify the weight assignment $\wt_{X,Y,z}$ for $L^{(1)}$, we choose the orthogonal coordinate system $Oij$ so that the horizontal $i$-axis runs along the base of the region, and the vertical $j$-axis passes through the lower vertex of the triangular component. Recall that a vertical lozenge at the horizontal coordinate $k$ is weighted by $f_{X,Y}(k+z)=\frac{Xq^{k+z}+Yq^{-k-z}}{2}$, and all lozenges of the other two orientations are weighted by $1$.

The region $\overline{L}^{(1)}_{n,x}(P, Q, R, S,C)$ is defined similarly. We start with the  parallelogram of side-lengths $n+x,n,n+x,n$.A diagonal
cut from the upper-right vertex to the base divides the parallelogram into an up-pointing triangle of
side length $n$  and a trapezoid. Next, we remove unit triangles along the sides of the triangle and the trapezoid at the positions
specified by the same sets $P,Q,R,S$ as shown in the right picture in Figure \ref{Lshape1}. The barriers are also placed along the cut at the same positions in $C$. We denote by $\overline{L}^{(1)}=\overline{L}^{(1)}_{n,x}(P, Q, R, S,C)$ the resulting region. The weight assignment $\wt_{X,Y,z}$ on $\overline{L}^{(1)}$ is determined by the coordinate system $Oij$ with the $i$-axis running along the base of the region and the $j$-axis passing through the midpoints of the bases of the trapezoidal component.

\begin{lem}\label{lem3}
Let $n$ and $x$ be nonnegative integers. Let $P, Q, R, S,C$ be subsets of $[n]$ satisfying $|P|+|Q|=|R|+|S|\leq n$, $1\notin P\cap Q$,  $1\notin R\cap S$, and $C\cap (Q\cup R)=\emptyset$. If $L^{(1)}_{n,x}(P, Q, R, S)$ admits a lozenge tiling, then
\begin{align*}
\frac{M_{X,Y,z}(\overline{L}^{(1)}_{n,x}(P, Q, R, S,C))}{\M_{X,Y,z}(L^{(1)}_{n,x}(P, Q, R, S,C))}&=\frac{1}{E_{x,z}(S+x)} \frac{D_x(P)D_x(Q)}{D_x(R)D_x(S)}.
\end{align*}
\end{lem}

\begin{proof} For a tiling $\tau$ of $L^{(1)}=L^{(1)}_{n,x}(P, Q, R, S,C)$, let $U=U_{\tau}$ be the set of positions of the lozenges crossing our cut. We call $U_\tau$ the \emph{crossing set} of $\tau$. Removing these lozenges preserves the tiling weight because each has weight $1$. After the lozenge removal, we obtain a tiling $\tau'$ of ${L}^{(1)}_{n,x}(P, Q\cup U, R\cup U, S,C)$. The tiling $\tau'$ is naturally partitioned into two tilings of the two component regions, denoted by $R_1$ and $R_2$.  The two components are congruent to the semihexagons described in the shuffling lemmas (Lemmas \ref{shufflingthm1} and \ref{shufflingthm2}) above. To identify their generating functions, we must also match their coordinate systems. We therefore apply the reflection and translation identities before using the shuffling lemmas. In particular, the trapezoidal  component $R_1$  is congruent to the region $S'_{n,x}(R\cup U, S)$  in Lemma \ref{shufflingthm2}. However, we need to translate it horizontally by $x$ units to the left to align the coordinate systems. By property $(P3)$ of the weight assignment $\wt_{X,Y,z}$, we have
 \[\M_{X,Y,z}(R_1)=\M_{X,Y,x+z}(S'_{n,x}(R\cup U, S)) .\]
The triangular component $R_2$ is obtained from reflecting the region $S_{n,0}(P, Q\cup U)$ over the $i$-axis, and the coordinate systems agree. By the property (P2) of the weight assignment, we have
\[\M_{X,Y,z}(R_2)=\M_{X,Y,z}(S_{n,0}(P,Q\cup U)).\]
We can therefore write
\begin{align}\label{Leq3}
\M_{X,Y,z}(L_{n,x}(P, Q\cup U, R\cup U, S,C))&= \M_{X,Y,z}(R_1) \cdot \M_{X,Y,z}(R_2)\\
&=\M_{X,Y,x+z}(S'_{n,x}(R\cup U, S)) \cdot \M_{X,Y,z}(S_{n,0}(P,Q\cup U)).
\end{align}

Let $U$ be a crossing set of a tiling of $L^{(1)}$. It is easy to see that the crossing lozenges must avoid the removed unit triangles and the barriers along the cut, and, by Lemma \ref{lem0},  the dent sets of both residual components must also be compatible. This implies that
\[U\subseteq [n]\setminus(Q\cup R\cup C)\quad \text{and} \quad |U|=n-|P|-|Q|=n-|R|-|S| \tag{$D_1$}\]
and
 \[\text{the pairs $(P,Q\cup U)$ and $(R\cup U,S)$ are compatible} \tag{$D_2$}\] 
If  $U$  satisfies the conditions ($D_1$) and ($D_2$), we say it is \emph{compatible} with the region $L^{(1)}$.

Conversely, suppose that $U$ satisfies these conditions, i.e., it is \emph{compatible} with the region $L^{(1)}$ . Each component then has a tiling, and adjoining the crossing lozenges at the positions in $U$ gives a tiling of the original region with crossing set $U$.
Thus the admissible crossing sets are exactly those described by these conditions.
Partitioning the tilings by their crossing sets, and using the fact that each crossing lozenge has
weight $1$, we obtain
\begin{align}\label{Leq1}
\M_{X,Y,z}(L^{(1)})=\sum_{U}\M_{X,Y,z}({L}^{(1)}_{n,x}(P, Q\cup U, R\cup U, S,C)),
\end{align}
where the sum is taken over all subsets $U$ of $[n]$ compatible with $L^{(1)}$.  

Similarly, for a  tiling $\overline{\tau}$ of $\overline{L}^{(1)}=\overline{L}^{(1)}_{n,x}(P, Q, R, S,C)$, define  $U=U_{\overline{\tau}}$ to be the crossing set. We also remove all crossing lozenges to get the tiling $\overline{\tau}'$ of the region  $\overline{L}^{(1)}_{n,x}(P, Q\cup U, R\cup U, S,C)$. The tiling  $\overline{\tau}'$ is partitioned into two tilings of two component regions: the triangle on the right is obtained from the region $S'_{n,0}(R\cup U, S)$ in Lemma \ref{shufflingthm2} by translating it by $x$ units along the $i$-axis (to match the coordinate); the trapezoid on the left is obtained from reflecting the region $S_{n,x}(P, Q\cup U)$ over the $i$-axis (and the coordinate systems already match). This partition gives
\begin{align}\label{Leq5}
\M_{X,Y,z}(\overline{L}_{n,x}^{(1)}(P, Q\cup U, R\cup U, S,C))=\M_{X,Y,x+z}(S'_{n,0}(R\cup U, S)) \cdot \M_{X,Y,z}(S_{n,x}(P,Q\cup U)).
\end{align}
 By Lemma \ref{lem0}, the set $U_{\overline{\tau}}$ must satisfy the same conditions ($D_1$) and ($D_2$). Conversely, for every set $U$ satisfying these conditions, there exists a tiling $\overline{\tau}$ whose crossing set is $U$.
Thus, we can write
\begin{align}\label{Leq2}
\M_{X,Y,z}(\overline{L}^{(1)})=\sum_{U}\M_{X,Y,z}(\overline{L}^{(1)}_{n,x}(P, Q\cup U, R\cup U, S,C)),
\end{align}
where the sum is taken over the same compatible position sets $U$ as in (\ref{Leq1}). 

  \begin{figure}[ht]\centering
   \includegraphics[width=15cm]{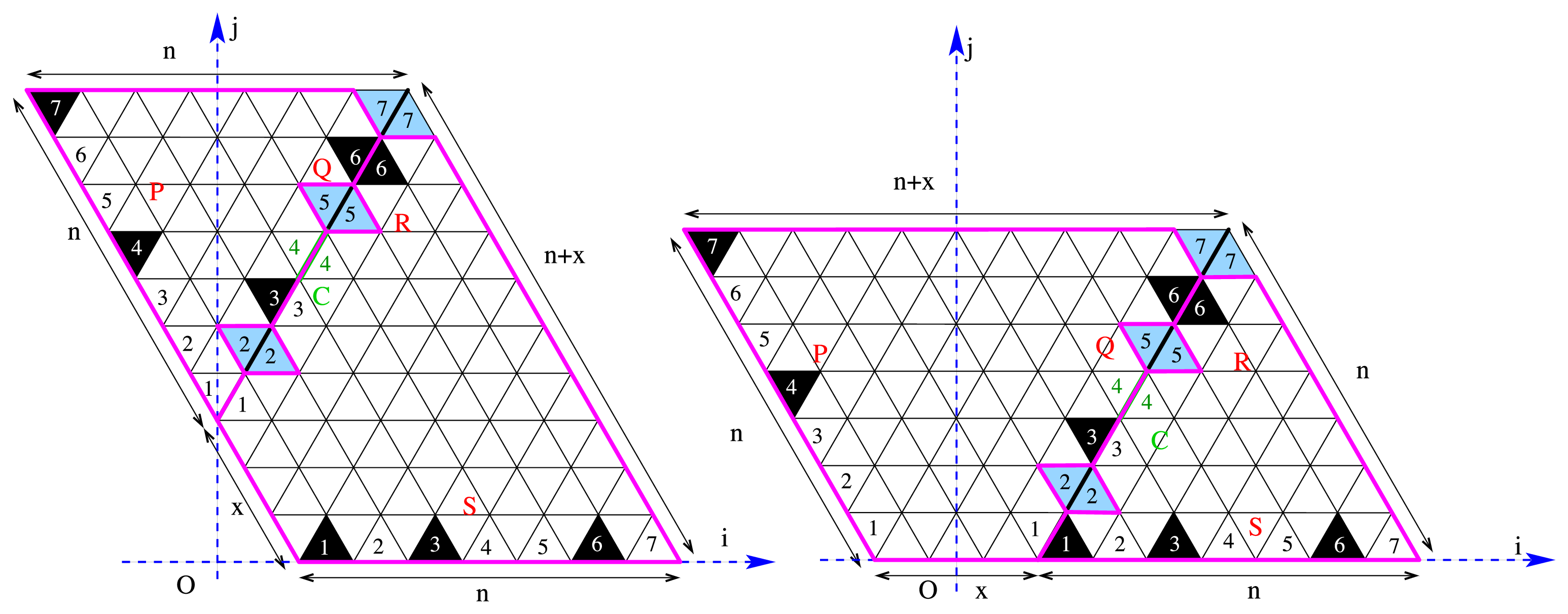}
 \caption{Decomposition of the regions in Figure \ref{Lshape1} after fixing the crossing set $U=\{2,5,7\}$. Removing the shaded crossing lozenges separates each region into two semihexagon components.} \label{Lshape1b}
 \end{figure}
 
By the Lemma \ref{shufflingthm1} and \ref{shufflingthm2} and equations (\ref{Leq3}) and (\ref{Leq5}), we have
\begin{align}\label{Leq6}
&\frac{\M_{X,Y,z}(\overline{L}_{n,x}(P, Q\cup U, R\cup U, S,C))}{\M_{X,Y,z}(L_{n,x}(P, Q\cup U, R\cup U, S,C))}=\frac{\M_{X,Y,x+z}(S'_{n,0}(R\cup U, S))}{\M_{X,Y,x+z}(S'_{n,x}(R\cup U, S))} \cdot \frac{\M_{X,Y,z}(S_{n,x}(P,Q\cup U))}{\M_{X,Y,z}(S_{n,o}(P,Q\cup U))}\notag\\
&=\frac{F_x(n)}{E_{x,z}(S+x)} \frac{1}{D_x(R \cup U)D_x(S)}\cdot \frac{1}{F_x(n)}D_x(Q\cup U)D_x(P)\notag\\
&=\frac{1}{E_{x,z}(S+x)}\frac{D_x(Q\cup U)D_x(P)}{D_x(R \cup U)D_x(S)}\notag\\
&=\frac{1}{E_{x,z}(S+x)}\frac{D_x(P)D_x(Q)}{D_x(R)D_x(S)}.
\end{align}
Here we use the following identities 
 \[D_x(Q\cup U)=D_x(Q) D_x(U)\quad \text{and} \quad D_x(R\cup U)=D_x(R) D_x(U)\] 
 (since $U$ is disjoint from $Q$ and $R$) to cancel the factors indexed by elements in $U$.
The resulting ratio is independent of $U$. 

Since the ratio is independent of $U$, multiplying each termwise ratio identity by its denominator, summing
the resulting equalities over $U$, and applying (\ref{Leq1}) and (\ref{Leq2}), we obtain
\begin{align}
&\frac{\M_{X,Y,z}(\overline{L}_{n,x}(P, Q, R, S,C))}{\M_{X,Y,z}(L_{n,x}(P, Q, R, S,C))} =\frac{1}{E_{x,z}(S+x)}\frac{D_x(P)D_x(Q)}{D_x(R)D_x(S)}.
\end{align}
\end{proof}

\begin{figure}[ht]\centering
    \includegraphics[width=15cm]{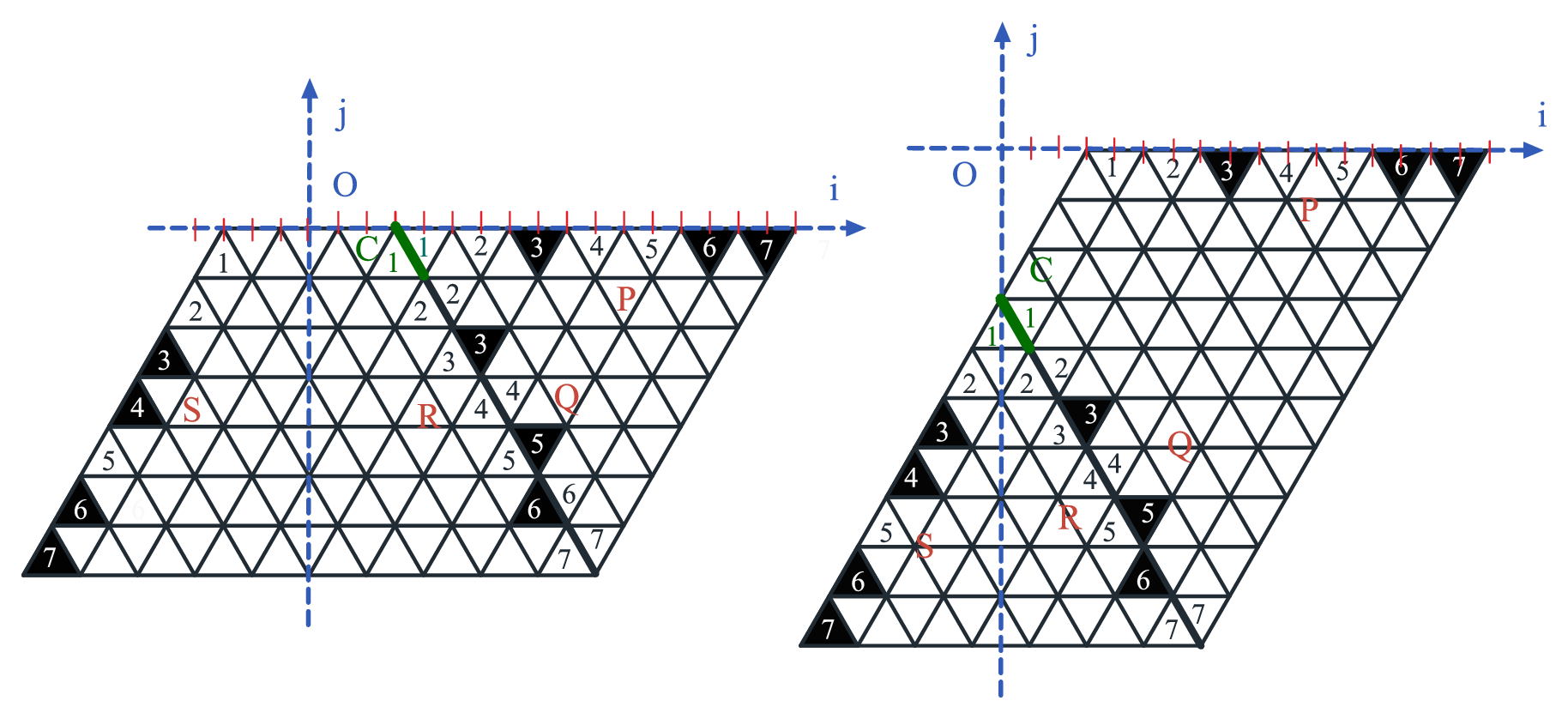}
 \caption{The regions $L^{(2)}_{7,3}(P,Q,R,S,C)$ (left) and  $\overline{L}^{(2)}_{7,3}(P,Q,R,S,C)$ (right) with $P=\{3,6,7\}$, $Q=\{3,5\}$, $R=\{6\}$, $S=\{3,4,6,7\}$, and $C=\{1\}$. The coordinate systems show the weight assignments used in
Lemma \ref{lem4}.} \label{Lshape2}
 \end{figure}

 Next, we consider a pair of similar parallelogram-shaped regions. We start with a rightward-slanting parallelogram $L^{(2)}_{n,x}$ of side-lengths $n+x,n,n+x,n$. We make a diagonally cut from the lower-right vertex to the top side to divide the parallelogram into two components: a down-pointing triangle of side-length $n$ on the right and a trapezoid on the left. We also remove unit triangles at the positions specified by $P,Q,R,S$ and place barriers at the positions specified by $C$  along the cut. See the left picture in Figure \ref{Lshape2}. These five subsets of $[n]$ satisfy the same conditions imposed in the definition of the $L^{(1)}$-type regions. We denote the region by $L^{(2)}_{n,x}(P, Q, R, S,C)$. For the weight assignment, we introduce the coordinate system $Oij$ so that the horizontal $i$-axis runs along the top of the region and the vertical $j$-axis passes through the midpoints of two parallel sides of the trapezoidal component.
 
 The region  $\overline{L}^{(2)}_{n,x}(P, Q, R, S,C)$ is defined similarly from the parallelogram $L^{(2)}_{n,x}$ of side-lengths $n,n+x,n,n+x$ with the diagonal cut running from the lower-right vertex to the oblique side on the left. This divides the region into a triangle and a trapezoid.  Remove unit triangles at the positions specified by the sets $P,Q,R,S$ on the sides of the component regions and place barriers at the positions specified by $C$ along the cut. See the right picture in Figure \ref{Lshape2}. Finally, we place  the coordinate system  $Oij$ so that the $i$-axis runs along the top side and the $j$-axis passes the top vertex of the triangular component.
 
\begin{lem}\label{lem4}
Let $n$ and $x$ be nonnegative integers. Suppose $P, Q, R, S,C$  are five subsets of $[n]$ satisfying the same hypotheses of Lemma \ref{lem3}. If $L^{(2)}_{n,x}(P, Q, R, S,C)$ admits a lozenge tiling, then
\begin{align*}
\frac{\M_{X,Y,z}(\overline{L}^{(2)}_{n,x}(P, Q, R, S,C))}{\M_{X,Y,z}(L^{(2)}_{n,x}(P, Q, R, S,C))}=&E_{x,z}(P+x)\cdot \frac{D_x(P)D_x(Q)}{D_x(R)D_x(S)}.
\end{align*}
\end{lem}

\begin{proof}
We now flip the regions $L^{(2)}_{n,x}(P, Q, R, S,C)$ and $\overline{L}^{(2)}_{n,x}(P, Q, R, S,C)$  over the horizontal $i$-axis to get the regions $\overline{L}^{(1)}_{n,x}(S, R, Q, P,C)$ and $L^{(1)}_{n,x}(S, R, Q, P,C)$ in Lemma \ref{lem3}, respectively. By property (P3), , reflection across the horizontal i-axis preserves the tiling generating functions, we have
\[\M_{X,Y,z}(L^{(2)}_{n,x}(P, Q, R, S,C))=\M_{X,Y,z}(\overline{L}^{(1)}_{n,x}(S, R, Q, P,C))\]
\[\M_{X,Y,z}(\overline{L}^{(2)}_{n,x}(P, Q, R, S,C))=\M_{X,Y,z}(L^{(1)}_{n,x}(S, R, Q, P,C)).\]
 Then the lemma now follows from Lemma \ref{lem3}.
\end{proof}

   \begin{figure}[ht]\centering
       \includegraphics[width=15cm]{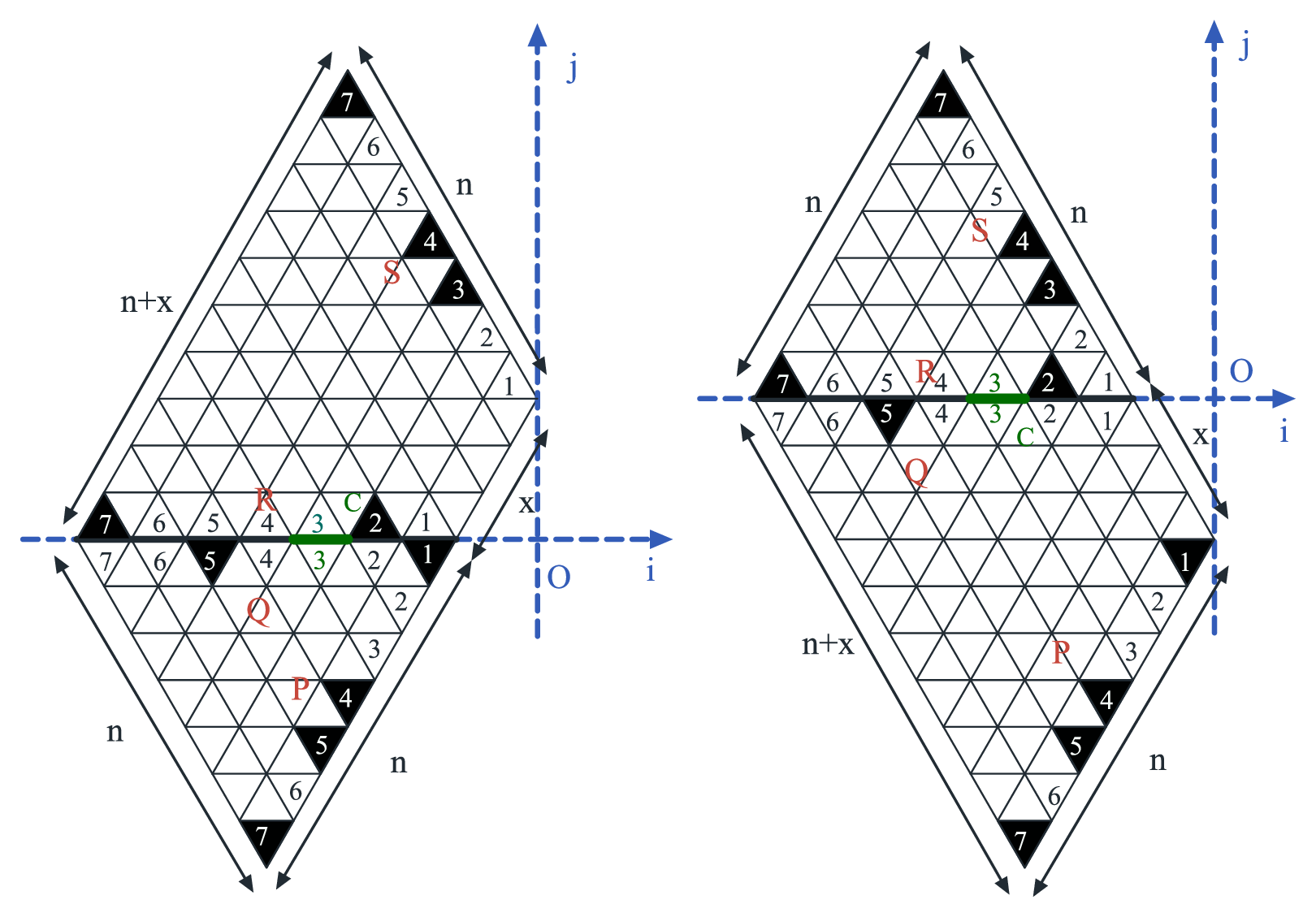}
 \caption{The regions $L^{(3)}_{7,3}(P,Q,R,S,C)$ (left) and  $\overline{L}^{(3)}_{7,3}(P,Q,R,S,C)$ (right) with $P=\{1,4,5,7\}$, $Q=\{5\}$, $R=\{2,7\}$, $S=\{3,4,7\}$, and $C=\{3\}$. Each region is divided by a horizontal cut into a triangle and a
trapezoid. The coordinate systems determine the position-dependent weights.} \label{Lshape3}
 \end{figure}

Our third pair of parallelogram-shaped regions has different orientations. We start with the parallelogram $L^{(3)}_{n,x}$ of side-lengths $n+x,n,n+x,n$ in clockwise order starting from the northwest side as shown in Figure \ref{Lshape3} (left). Our region also has a \emph{horizontal} cut from the left vertex to a point on the southeast side. The cut also divides the region into two components: a trapezoid on top and a down-pointing triangle of side-length $n$ below it. We also remove the unit triangles along the sides of the two components at the position specified by $P,$ $Q,$ $R,$ and $S$ and place barriers at the positions specified by $C$ along the cut. The sets $P,Q,R,S,C$  are the subsets of $[n]$ satisfying the same requirements as in the definition of the first pair of the  $L^{(1)}$- and $\overline{L}^{(1)}$-regions. We denote the resulting region by $L^{(3)}_{n,x}(P, Q, R, S,C)$. Next, we introduce the coordinate system $Oij$ so that the $i$-axis runs horizontally along the cut and the $j$-axis passes through the right vertex of the region. The lozenges are weighted by the weight assignment $\wt_{X,Y,z}$.

The region $\overline{L}^{(3)}_{n,x}(P, Q, R, S,C)$ is defined similarly from the parallelogram  $\overline{L}^{(3)}_{n,x}$ of side-lengths $n,n+x,n,n+x$   in clockwise order starting from the northwest side as shown in Figure \ref{Lshape3} (right). We make the cut and remove the unit triangles at the positions specified by the sets $P,Q,R,S$ and place barriers at the positions specified by the set $C$ along the cut.  The resulting region is denoted by $\overline{L}^{(3)}_{n,x}(P, Q, R, S,C)$. The coordinate system $OIj$ is positioned in the same way. We use the same weight assignment $\wt_{X,Y,z}$.
 
\begin{lem}\label{lem5}
Let $n$ and $x$ be nonnegative integers. Suppose $P, Q, R, S,C$ are five subsets of $[n]$ satisfying the same hypotheses of Lemma \ref{lem3}. If $L^{(3)}_{n,x}(P, Q, R, S,C)$ admits a lozenge tiling, then
\begin{align*}
\frac{\M_{X,Y,z}(\overline{L}^{(3)}_{n,x}(P, Q, R, S,C))}{\M_{X,Y,z}(L^{(3)}_{n,x}(P, Q, R, S,C))}=&\frac{E'_{x,z}(Q+x)}{E'_{x,z}(R+x)} \frac{D_x(P)D_x(Q)}{D_x(R)D_x(S)}.
\end{align*}
\end{lem}

   \begin{figure}[ht]\centering
         \includegraphics[width=15cm]{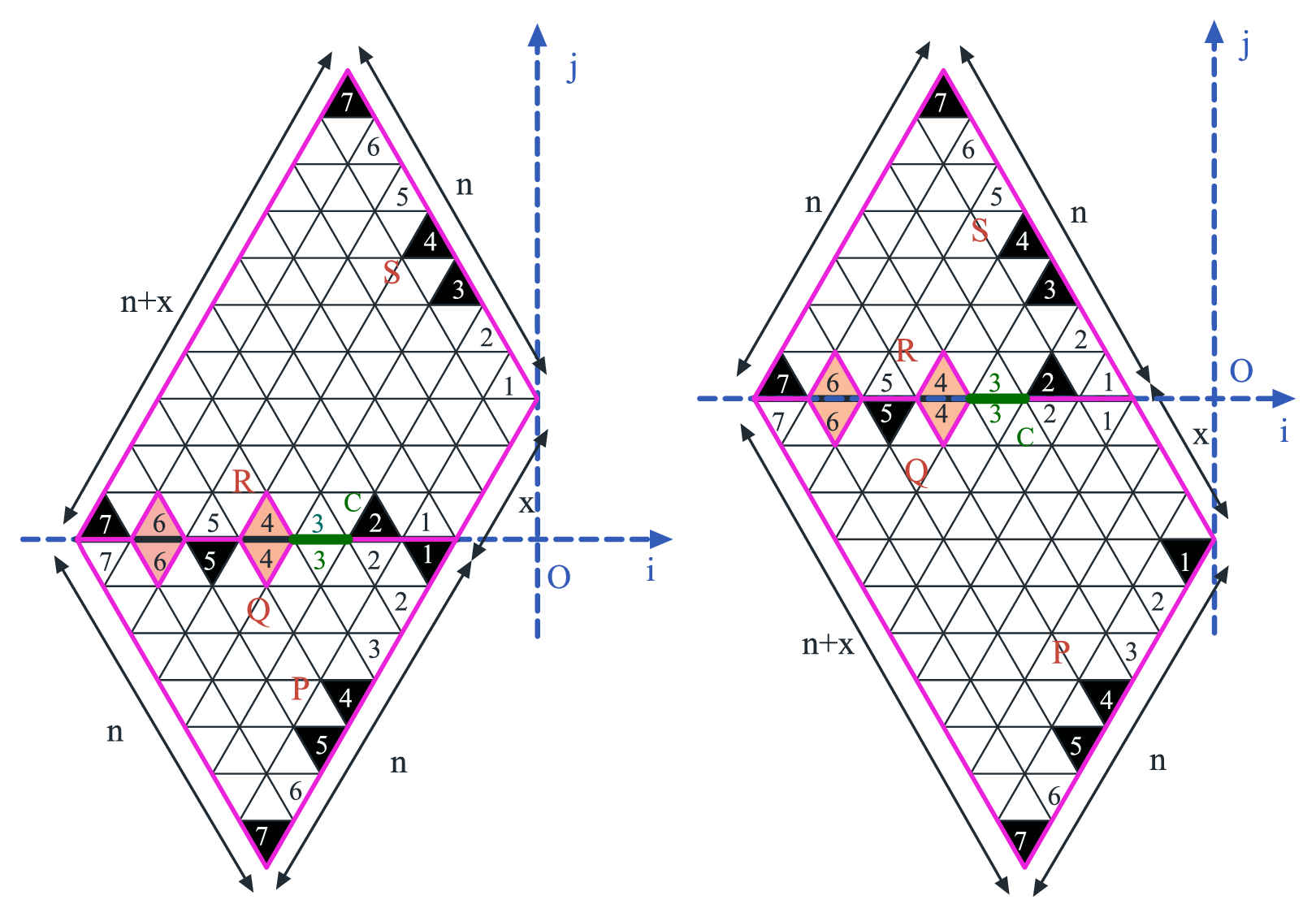}
 \caption{Decompositions of the regions in Figure \ref{Lshape3} and with the crossing set $U=\{4,6\}$.  The shaded vertical lozenges cross the horizontal cut. Their weights must be included when the generating functions are decomposed according to $U$.} \label{Lshape3b}
 \end{figure}

\begin{proof} The proof is similar to that of  Lemma \ref{lem3}.

For a tiling of $L^{(3)}_{n,x}=L^{(3}_{n,x}(P,Q,R,S,C)$, let $U$ be the set of positions of the lozenges crossing the horizontal cut. Similarly, $U$ must satisfy the same conditions ($D_1$) and ($D_2$). We also call such a set $U$ compatible. Moreover, for each compatible set $U$, there is a tiling of $L^{(3)}_{n,x}$ whose crossing set is $U$.
As in the proof of Lemma \ref{lem3}, partitioning the tilings by their crossing sets gives
\begin{align}\label{Leq4b}
\M_{X,Y,z}(L^{(3)})=\sum_{U}\M_{X,Y,z}({L}^{(3)}_{n,x}(P, Q\cup U, R\cup U, S,C)) \cdot \wt(U)
\end{align}
where the sum is taken over all compatible sets $U$ and $\wt(U)$ is the product of the weights of the crossing lozenges at the positions specified by $U$. This case is different from the case of $L^{(1)}$-regions in Lemma \ref{lem3}, as the crossing lozenges are vertical lozenges and their weights are not identically equal to $1$.
Similarly, we write
\begin{align}\label{Leq5b}
\M_{X,Y,z}(\overline{L}^{(3)})=\sum_{U}\M_{X,Y,z}(\overline{L}^{(3)}_{n,x}(P, Q\cup U, R\cup U, S,C)) \cdot \overline{\wt}(U)
\end{align}
where the sum is over the same collection of compatible sets $U$ as in (\ref{Leq4b}) and $\overline{\wt}(U)$ is the product of weight of all crossing lozenges at the positions specified by $U$. Investigate the $i$-coordinates of the crossing lozenges specified by $U$ in the two regions $L^{(3)}$ and $\overline{L}^{(3)}$, we actually have 
\[\wt(U)=\overline{\wt}(U)=\prod_{u\in U} f_{X,Y}(z-x-2u+1).\]

For a given compatible set $U$,  we divide ${L}^{(3)}_{n,x}(P, Q\cup U, R\cup U, S,C)$ along the cut into two  components: a triangle of side length $n$ and a trapezoid, denoted by $R_1$ and $R_2$, respectively. 
By the definition of $U$, each tiling of ${L}^{(3)}_{n,x}(P, Q\cup U, R\cup U, S,C)$ can be partitioned into two tilings of the triangle $R_1$ and  the trapezoid $R_2$. We then have
\begin{align}\label{Leq5d}
 \M_{X,Y,z}({L}^{(3)}_{n,x}(P, Q\cup U, R\cup U, S,C))= \M_{X,Y,z}(R_1) \cdot  \M_{X,Y,z}(R_2).
\end{align}

We reflect the triangle $R_1$ is obtained by translating horizontally the region $S'_{n,0}(P,Q\cup U)$ Lemma  \ref{shufflingthm2} by $x$ units to the right, reflecting over the $i$-axis, and then reflecting over the $j$-axis. Denote by $R'$ and $R''$ the resulting regions after the first and the second transformation above. It means $R_1$ is obtained by reflecting $R''$ over the $j$-axis, $R''$ is obtained by reflecting $R'$ over the $i$-axis, and $R'$ is obtained by translating the region $S'_{n,0}(P,Q\cup U)$ by $x$ units to the right. By the property (P1), (P2), and (P3), respectively, we have
\begin{align*}
\M_{X,Y,z}(R_1)=\M_{Y,X,-z}(R''),
\end{align*}
\begin{align*}
\M_{Y,X,-z}(R'')=\M_{Y,X,-z}(R'),
\end{align*}
and
\begin{align*}
\M_{Y,X,-z}(R')=\M_{Y,X,x-z}(R_1).
\end{align*}
The three identities imply that
\begin{align}\label{Leq5e}
\M_{X,Y,z}(R_1)=\M_{Y,X,x-z}(S'_{n,0}(P,Q\cup U)).
\end{align}

The  trapezoid $R_2$ is obtained by translating $S'_{n,x}(S, R\cup U)$ horizontally $x$ units to the right and reflecting over the $j$-axis. If we call the intermediate region is $\widehat{R}$, the by property (P1), we have 
\begin{align*}
\M_{X,Y,z}(R_2)=\M_{Y,X,-z}(\widehat{R}). 
\end{align*}
By property (P3), we have
\begin{align*}
\M_{Y,X,-z}(\widehat{R})=\M_{Y,X,x-z}(S'_{n,x}(S, R\cup U)). 
\end{align*}
Thus, 
\begin{align}\label{Leq5f}
\M_{X,Y,z}(R_2)=\M_{Y,X,x-z}(S'_{n,x}(S, R\cup U)). 
\end{align}
By (\ref{Leq5d}), (\ref{Leq5e}), and (\ref{Leq5d}), we obtain
\begin{align*}
\M_{X,Y,z}(L^{(3)}_{n,x}(P, Q\cup U, R\cup U, S,C))=\M_{Y,X,x-z}(S'_{n,0}( P, Q\cup U)) \cdot \M_{Y,X,x-z}(S'_{n,x}(S,R\cup U)).
\end{align*}
Similarly,  we also divide $\overline{L}^{(3)}_{n,x}(P, Q\cup U, R\cup U, S,C)$ along the cut into two components: a triangle and a trapezoid. The trapezoid can be obtained from   $S'_{n,x}(P,Q\cup U)$ in Lemma  \ref{shufflingthm2} by translating it by $x$ units along $i$-axis, reflecting it over $i$-axis, and reflecting it over the $j$-axis. The triangle can be obtained from $S'_{n,0}(S, R\cup U)$ by translating by $x$ units along $i$-axis and reflecting it over the $j$-axis. We can also write
\begin{align*}
\M_{X,Y,z}(\overline{L}^{(3)}_{n,x}(P, Q\cup U, R\cup U, S,C))=\M_{Y,X,x-z}(S'_{n,x}( P, Q\cup U)) \cdot \M_{Y,X,x-z}(S'_{n,0}(S,R\cup U)).
\end{align*}

For each compatible set $U$, Lemma \ref{shufflingthm2} gives
\begin{align}\label{Leq8}
&\frac{\M_{X,Y,z}(\overline{L}^{(3)}_{n,x}(P, Q\cup U, R\cup U, S,C))}{\M_{X,Y,z}(L^{(3)}_{n,x}(P, Q\cup U, R\cup U , S,C))}=\frac{\M_{Y,X,x-z}(S'_{n,x}( P, Q\cup U))}{\M_{Y,X,x-z}(S'_{n,0}( P, Q\cup U))}\cdot \frac{\M_{Y,X,x-z}(S'_{n,0}(S,R\cup U))}{ \M_{Y,X,x-z}(S'_{n,x}(S,R\cup U))}\notag\\
&=\frac{E'_{x,z}((Q\cup U) +x)}{F_x(n)}D_x(P)D_x(Q\cup U) \cdot \frac{F_x(n)}{E'_{x,z}((R\cup U)+x)} \frac{1}{D_x(R\cup U)D_x(S)} \notag\\
&=\frac{E'_{x,z}((Q\cup U) +x)}{E'_{x,z}((R\cup U)+x)}\frac{D_x(P)D_x(Q\cup U)}{D_x(R\cup U)D_x(S)}\notag\\
&=\frac{E'_{x,z}(Q +x)}{E'_{x,z}(R+x)}\frac{D_x(P)D_x(Q)}{D_x(R)D_x(S)}.
\end{align}
Here we use the identities
\[D_x(Q\cup U)=D_x(Q) D_x(U), \quad D_x(R\cup U)=D_x(R) D_x(U),\] 
\[\frac{E'_{x,z}((Q\cup U)+x)}{E'_{x,z}((R\cup U)+x)}=\frac{E'_{x,z}(Q+x)E'_{x,z}(U+x)}{E'_{x,z}(R+x)E'_{x,z}(U+x)}=\frac{E'_{x,z}(Q+x)}{E'_{x,z}(R+x)}\]  (since $U$ is disjoint from $Q$ and $R$)  to cancel factors indexed elements in the set $U$. 
The resulting ratio is independent of $U$. 

Since the ratio is independent of $U$, multiplying each termwise ratio identity by its denominator, multiplying both sides by the factor $\wt(U)=\overline{\wt}(U)$, summing
the resulting equalities over $U$, and applying (\ref{Leq4b}) and (\ref{Leq5b}), we have
\begin{align}
\frac{\M_{X,Y,z}(\overline{L}^{(3)}_{n,x}(P, Q, R, S,C))}{\M_{X,Y,z}(L^{(3)}_{n,x}(P, Q, R, S,C))}=\frac{E'_{x,z}(Q +x)}{E'_{x,z}(R+x)}\frac{D_x(P)D_x(Q)}{D_x(R)D_x(S)}.
\end{align}
\end{proof}

\begin{cor}\label{cor}
Assume that $n$ and $x$ are two nonnegative integers. Assume that $P,Q,R,S,C$ are five subsets of $[n]$ satisfying the same hypotheses as in Lemma \ref{lem3}. Then the region $L^{(i)}(P, Q, R, S,C)$ admits a tiling if and only if the region $\overline{L}^{(i)}(P, Q, R, S,C)$ admits a tiling, for $i=1,2,3$.
\end{cor}
\begin{proof}
The case $i=1$ follows directly from the identities (\ref{Leq1}) and (\ref{Leq2}). The case $i=2$ follows from the case $i=1$ by reflection over the $i$-axis. The case $i=3$ follows directly from the identities (\ref{Leq4b}) and (\ref{Leq5b}).
\end{proof}

We now apply the component identities in Lemmas \ref{lem3},  \ref{lem4},  and \ref{lem5} to prove our main theorem.

   \begin{figure}[ht]\centering
         \includegraphics[width=15cm]{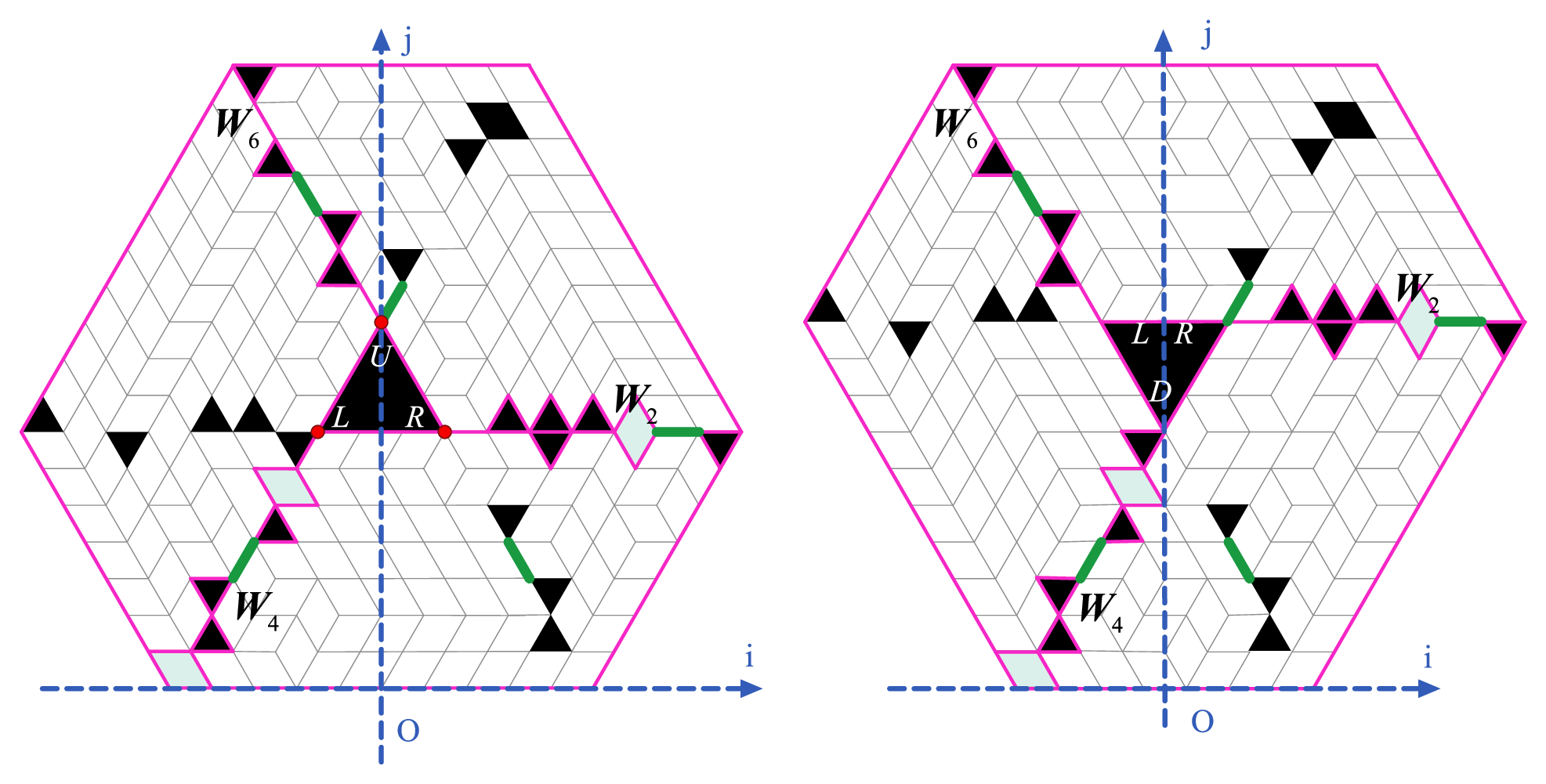}
 \caption{(a) Decomposing a residual tiling of  $H_{7,3}(\mathbf{A},\mathbf{B},\mathbf{C})$ contributing to $Z_{\mathbf{W}}$ into three regions of  types $L^{(1)}(\mathbf{W})$, $L^{(2)}(\mathbf{W})$, and $L^{(3)}(\mathbf{W})$ with $\mathbf{W}=(W_2,W_4,W_6)=(\{5\},\{2,7\},\emptyset)$. (b) Decomposing a residual tiling of $\overline{H}_{7,3}(\mathbf{A},\mathbf{B},\mathbf{C})$ contributing to $\overline{Z}_{\mathbf{W}}$ into three regions of  types $\overline{L}^{(1)}$, $\overline{L}^{(2)}$, and $\overline{L}^{(3)}$ with the same set $\mathbf{W}$. The shaded lozenges occupy the crossing positions recorded by $\mathbf{W}$.} \label{proofpic}
 \end{figure}

\begin{proof}[Proof of Theorem \ref{mainthm}]
We follow the decomposition used in Byun’s proof of \cite[Theorem 2.1]{Byun2}. We partition the  set of lozenge tilings of the snowflake region $H_{n,x}(\mathbf{A},\mathbf{B},\mathbf{C})$ according to the positions of lozenges crossing the axes joining the vertices of the central triangular hole to the upper-left, right, and lower-left vertices of the outer hexagon. For a fixed tiling $\pi$ of $H_{n,x}(\mathbf{A},\mathbf{B},\mathbf{C})$,  let $W_2,W_4,W_6$  record the positions of  all lozenges of $\pi$ crossing the three selected axes.  

Since the crossing lozenges must avoid the removed unit triangles and the barriers on the same axis, we must have
\begin{equation}\label{finalcondition}
W_{i}\cap (C_{i+1} \cup  A_{i+1} \cup B_{i})=\emptyset,
\end{equation}
for $i=2,4,6$, and $A_{7}=A_1$ and $C_7=C_1$ by convention. We assume this condition throughout. We call $\mathbf{W}=\mathbf{W_\pi}:=(W_2,W_4,W_6)$ the \emph{crossing set} of $\pi$.
Denote by  $c(\mathbf{W})$ the product of the  weights of the crossing lozenges determined by $\mathbf{W}$  in the region $H_{n,x}(\mathbf{A},\mathbf{B},\mathbf{C})$.
We denote by $Z_{\mathbf{W}}$ the weighted sum of all residual tilings after removing the lozenges  specified by $\mathbf{W}$ and
\emph{forbidding further crossings of the three selected axes}. Thus\footnote{Denote by $H_{n,x}(\mathbf{A}_\mathbf{W},\mathbf{B}_\mathbf{W},\mathbf{C})$ the residual region obtained from $H_{n,x}(\mathbf{A},\mathbf{B},\mathbf{C})$ by removing the crossing lozenges specified by $\mathbf{W}$. We \emph{cannot} simply write \begin{align*}
\M_{X,Y,z}(H_{n,x}(\mathbf{A},\mathbf{B},\mathbf{C}))=\sum_{\mathbf{W}} \M_{X,Y,z}(H_{n,x}(\mathbf{A}_\mathbf{W},\mathbf{B}^\mathbf{W},\mathbf{C})) \cdot c(\mathbf{W}),
\end{align*} where the sum is over all triples $\mathbf{W}$ satisfying  the condition (\ref{finalcondition}) and  for which the residual region $H_{n,x}(\mathbf{A}_\mathbf{W},\mathbf{B}_\mathbf{W},\mathbf{C})$ admits a tiling. The reason is that, the region $H_{n,x}(\mathbf{A}_\mathbf{W},\mathbf{B}_\mathbf{W},\mathbf{C})$ may have a tiling $\mu$ with further crossing lozenges, then the corresponding tiling of the original snowflake region is counted more than once in the sum on the right-hand side. Unlike the crossing set of an $L^{(i)}$-region, the sizes of the components $W_2,W_4,W_6$ of our crossing set $\mathbf{W}$ are not determined uniquely in a given snowflake region.}
\begin{align}\label{finaleq1}
\M_{X,Y,z}(H_{n,x}(\mathbf{A},\mathbf{B},\mathbf{C}))=\sum_{\mathbf{W}} Z_{\mathbf{W}} \cdot c(\mathbf{W}).
\end{align}


We call the triple  $\mathbf{W}=(W_2,W_4,W_6)$ \emph{compatible} with $H_{n,x}(\mathbf{A},\mathbf{B},\mathbf{C})$ if there is a tiling of the region whose crossing set is $\mathbf{W}$. This way, the sum in \ref{finaleq1} is taken over all the set of compatible triples $\mathbf{W}$.

For a fixed compatible $\mathbf{W}$,  each residual tiling contributing to $Z_{\mathbf{W}}$ is partitioned into three tilings of three parallelogram-shaped regions $L^{(1)}_{n,x}(\mathbf{W})=L^{(1)}_{n,x}(A_1\cup W_6, B_1,A_2,B_2\cup W_2,C_2)$,  $L^{(2)}_{n,x}(\mathbf{W})=L^{(2)}_{n,x}(A_3\cup W_2, B_3,A_4,B_4\cup W_4,C_4)$, and $L^{(3)}_{n,x}(\mathbf{W})=L^{(3)}_{n,x}(A_5\cup W_4, B_5,A_6,B_6\cup W_6,C_6)$.  See the left picture in Figure \ref{proofpic} for an example; the component regions are the components enclosed by the purple contours. 
We have 
\begin{align}\label{finaleq1b}
Z_{\mathbf{W}}=\M_{X,Y,z}(L^{(1)}_{n,x}(\mathbf{W})) \cdot  \M_{X,Y,z}(L^{(2)}_{n,x}(\mathbf{W})) \cdot \M_{X,Y,z}(L^{(3)}_{n,x}(\mathbf{W})).
\end{align}
Thus,
\begin{align}\label{finaleq1c}
\M_{X,Y,z}(H_{n,x}(\mathbf{A},\mathbf{B},\mathbf{C}))=\sum_{\mathbf{W}} \M_{X,Y,z}(L^{(1)}_{n,x}(\mathbf{W})) \cdot  \M_{X,Y,z}(L^{(2)}_{n,x}(\mathbf{W})) \cdot \M_{X,Y,z}(L^{(3)}_{n,x}(\mathbf{W})) \cdot c(\mathbf{W}).
\end{align}

We define similarly the crossing set $\mathbf{W}_{\overline{\pi}}$ of a tiling $\overline{\pi}$ of $\overline{H}_{n,x}(\mathbf{A},\mathbf{B},\mathbf{C})$. We also say that $\mathbf{W}$ is compatible with $\overline{H}_{n,x}(\mathbf{A},\mathbf{B},\mathbf{C})$ if there is a tiling of the region whose crossing set is $\mathbf{W}$. We define $\overline{Z}(\mathbf{W})$ and $\overline{c}(\mathbf{W})$ similarly. 

We have
 \begin{align}\label{finaleq1x}
\M_{X,Y,z}(\overline{H}_{n,x}(\mathbf{A},\mathbf{B},\mathbf{C}))=\sum_{\mathbf{W}} \overline{Z}_{\mathbf{W}} \cdot \overline{c}(\mathbf{W})
\end{align}
and
\begin{align}\label{finaleq1y}
\overline{Z}_{\mathbf{W}}=\M_{X,Y,z}(\overline{L}^{(1)}_{n,x}(\mathbf{W})) \cdot  \M_{X,Y,z}(\overline{L}^{(2)}_{n,x}(\mathbf{W})) \cdot \M_{X,Y,z}(\overline{L}^{(3)}_{n,x}(\mathbf{W})),
\end{align}
where $\overline{L}^{(1)}(\mathbf{W})=\overline{L}^{(1)}_{n,x}(A_1\cup W_6, B_1,A_2,B_2\cup W_2,C_2)$,  $\overline{L}^{(2)}(\mathbf{W})=\overline{L}^{(2)}_{n,x}(A_3\cup W_2, B_3,A_4,B_4\cup W_4,C_4)$, and $\overline{L}^{(3)}(\mathbf{W})=\overline{L}^{(3)}_{n,x}(A_5\cup W_4, B_5,A_6,B_6\cup W_6,C_6)$. See the right picture in Figure \ref{proofpic} for an example; the component regions are the ones restricted by the purple contours.
Similarly, we can write 
\begin{align}\label{finaleq2c}
\M_{X,Y,z}(\overline{H}_{n,x}(\mathbf{A},\mathbf{B},\mathbf{C}))=\sum_{\mathbf{W}}\M_{X,Y,z}(\overline{L}^{(1)}_{n,x}(\mathbf{W})) \cdot  \M_{X,Y,z}(\overline{L}^{(2)}_{n,x}(\mathbf{W})) \cdot \M_{X,Y,z}(\overline{L}^{(3)}_{n,x}(\mathbf{W}) ) \cdot \overline{c}(\mathbf{W}),
\end{align}
where the sum is taken over all compatible sets $\mathbf{W}$ of $\overline{H}_{n,x}(\mathbf{A},\mathbf{B},\mathbf{C})$.

\bigskip

We next show out that the sets of compatible $\mathbf{W}$ of $H_{n,x}(\mathbf{A},\mathbf{B},\mathbf{C})$ and $\overline{H}_{n,x}(\mathbf{A},\mathbf{B},\mathbf{C})$ are the same. 
Indeed, a compatible triple $\mathbf{W}$ of $H_{n,x}(\mathbf{A},\mathbf{B},\mathbf{C})$ gives three tileable components $L^{(i)}_{n,x}(\mathbf{W})$, for $i=1,2,3$; Corollary 4.7 gives tileability of the sibling components $\overline{L}^{(i)}_{n,x}(\mathbf{W})$, for $i=1,2,3$; gluing those tilings with the prescribed crossings gives a tiling of the sibling snowflake $\overline{H}_{n,x}(\mathbf{A},\mathbf{B},\mathbf{C})$, so $\mathbf{W}$ is a compatible triple of $\overline{H}_{n,x}(\mathbf{A},\mathbf{B},\mathbf{C})$. Reverse the argument for the converse.

 \bigskip
 
 One readily sees that, for a fixed compatible triple $\mathbf{W}=(W_2,W_4,W_6)$, the products $c(\mathbf{W})$ and  $\overline{c}(\mathbf{W})$ depend on only the component set $W_2$, as the  crossing lozenges specified by $W_4$ and $W_6$ all have weight $1$. Investigate the $i$-coordinates of the crossing lozenges specified by $W_2$ in the two regions $H_{n,x}(\mathbf{A},\mathbf{B},\mathbf{C})$ and $\overline{H}_{n,x}(\mathbf{A},\mathbf{B},\mathbf{C})$, we actually have 
\[c(\mathbf{W})=\overline{c}(\mathbf{W})=\prod_{w\in W_2} f_{X,Y}(z+x+2w-1).\]

\bigskip



We note that the coordinate system in $H_{n,x}(\mathbf{A},\mathbf{B},\mathbf{C})$  agrees with the coordinates of $L^{(1)}_{n,x}(\mathbf{W})$, $L^{(2)}_{n,x}(\mathbf{W})$, and $L^{(3)}_{n,x}(\mathbf{W})$  in Lemmas \ref{lem3}, \ref{lem4}, and \ref{lem5}, respectively. Similarly, the coordinate system in $\overline{H}_{n,x}(\mathbf{A},\mathbf{B},\mathbf{C})$  agrees with the coordinates of $\overline{L}^{(1)}_{n,x}(\mathbf{W})$, $\overline{L}^{(2)}_{n,x}(\mathbf{W})$, and $\overline{L}^{(3)}_{n,x}(\mathbf{W})$  in Lemmas \ref{lem3}, \ref{lem4}, and \ref{lem5}, respectively.
 Applying Lemmas \ref{lem3}, \ref{lem4}, and \ref{lem5} gives
\begin{align}
\frac{\M_{X,Y,z}(\overline{L}^{(1)}(\mathbf{W}))}{\M_{X,Y,z}(L^{(1)}(\mathbf{W}))} =&\frac{1}{E_{x,z}((B_2\cup W_2)+x)}\frac{D_x(A_1\cup W_6)D_x(B_1)}{D_x(A_2)D_x(B_2\cup W_2)},
\end{align}
\begin{align}
\frac{\M_{X,Y,z}(\overline{L}^{(2)}(\mathbf{W}))}{\M_{X,Y,z}(L^{(2)}(\mathbf{W}))} =&E_{x,z}((A_3\cup W_2)+x) \frac{D_x(A_3\cup W_2)D_x(B_3)}{ D_x(A_4)D_x(B_4\cup W_4)},
\end{align}
\begin{align}
\frac{\M_{X,Y,z}(\overline{L}^{(3)}(\mathbf{W}))}{\M_{X,Y,z}(L^{(3)}(\mathbf{W}))} =&\frac{E'_{x,z}(B_5+x)}{E'_{x,z}(A_6+x)} \frac{D_x(A_5\cup W_4)D_x(B_5)}{D_x(A_6)D_x(B_6\cup W_6)}.
\end{align}
We will use the cyclic index conventions $A_{7}:=A_{1}$ and  $B_{7}:=B_1$. Multiplying the three identity and canceling out the common factor gives
\begin{align}\label{finaleq5}
&\frac{\M_{X,Y,z}(\overline{L}^{(1)}(\mathbf{W}))}{\M_{X,Y,z}(L^{(1)}(\mathbf{W}))} \frac{\M_{X,Y,z}(\overline{L}^{(2)}(\mathbf{W}))}{\M_{X,Y,z}(L^{(2)}(\mathbf{W}))}  \frac{\M_{X,Y,z}(\overline{L}^{(3)}(\mathbf{W}))}{\M_{X,Y,z}(L^{(3)}(\mathbf{W}))}\notag\\
&=\frac{E_{x,z}((A_3\cup W_2)+x)E'_{x,z}(B_5+x)}{E_{x,z}((B_2\cup W_2)+x)E'_{x,z}(A_6+x)} \prod_{j=1}^{3}\frac{D_x(A_{2j+1}\cup W_{2j})D_x(B_{2j+1})}{D_x(A_{2j})D_x(B_{2j}\cup W_{2j})}\notag\\
&=\frac{E_{x,z}(A_3+x)E'_{x,z}(B_5+x)}{E_{x,z}(B_2+x)E'_{x,z}(A_6+x)} \prod_{j=1}^{3}\frac{D_x(A_{2j+1})D_x(B_{2j+1})}{D_x(A_{2j})D_x(B_{2j})}\notag\\
&=\frac{E_{x,z}(A_3+x)E'_{x,z}(B_5+x)}{E_{x,z}(B_2+x)E'_{x,z}(A_6+x)} \prod_{j=1}^{3}\frac{D_x(A_{2j-1})D_x(B_{2j-1})}{D_x(A_{2j})D_x(B_{2j})}.
\end{align}
Here we use the identities
\[D_x(A_{2j+1}\cup W_{2j})=D_x(A_{2j+1})D_x(W_{2j}),  D_x(B_{2j}\cup W_{2j})=D_x(B_{2j})D_x(W_{2j}),\]
since $W_{2j}\cap A_{2j+1}=W_{2j}\cap B_{2j}=\emptyset$  to cancel the factors indexed by elements in $W_{2j}$.

Let $K$ be the expression on the right of (\ref{finaleq5}), independent of $\mathbf{W}$. Since $c(\mathbf{W})=\overline{c}(\mathbf{W})$,  by (\ref{finaleq1b}) and (\ref{finaleq1y}), we have
\[\overline{Z}_{\mathbf{W}}  \overline{c}(\mathbf{W})= K \cdot Z_{\mathbf{W} }c(\mathbf{W}),\]
for any common compatible corssing triple $\mathbf{W}$.
Summing over all such triples $\mathbf{W}$ and applying (\ref{finaleq1c}) and (\ref{finaleq2c}), we have
\begin{align}
\frac{\M_{X,Y,z}(\overline{H}_{n,x}(\mathbf{A},\mathbf{B},\mathbf{C}))}{\M_{X,Y,z}(H_{n,x}(\mathbf{A},\mathbf{B},\mathbf{C}))}&=\frac{E_{x,z}(A_3+x)E'_{x,z}(B_5+x)}{E_{x,z}(B_2+x)E'_{x,z}(A_6+x)}  \prod_{j=1}^{3}\frac{D_x(A_{2j-1})D_x(B_{2j-1})}{D_x(A_{2j})D_x(B_{2j})}.
\end{align}
This is (\ref{maineq}), completing the proof of Theorem \ref{mainthm}.
\end{proof}

\section{Proof of the Shuffling Lemmas}
In this section we prove Lemma \ref{shufflingthm2}. The proof of Lemma \ref{shufflingthm1} follows by adapting the argument in \cite{Semitwodent}
or the determinant evaluation in \cite{FulmekShuffling2}. We first record two weighted enumeration formulas for one-sided
dented semihexagons and derive the boundary identity in Lemma \ref{extralem}. We then apply Kuo condensation \cite{Kuo}
and induction to prove the desired multiplicative identity

 \begin{figure}[ht]\centering
          \includegraphics[width=15cm]{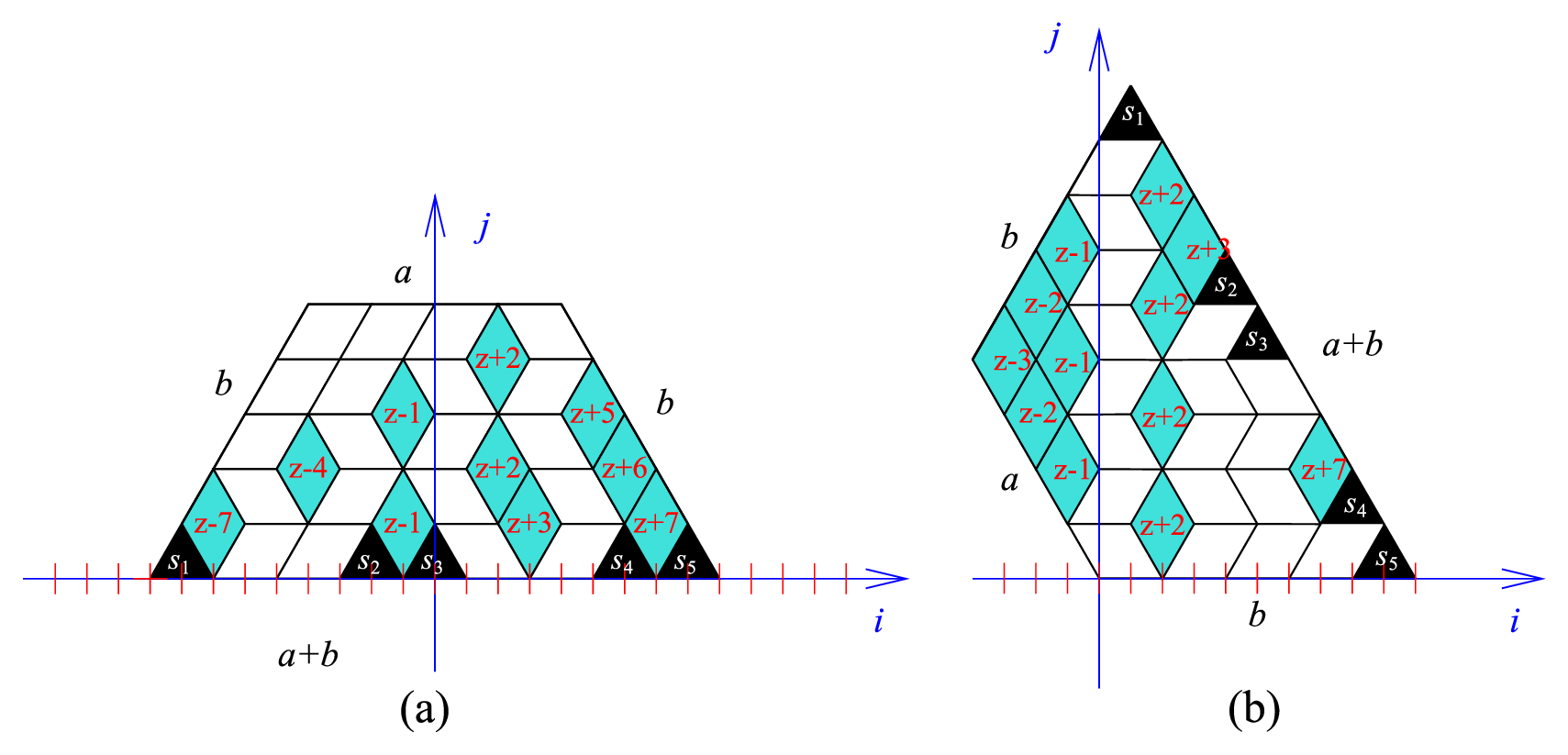}
 \caption{The semihexagon with removed triangles on the bases $SH_{a,b}(s_1,s_2,\dots,s_b)$ and $SH'_{a,b}(s_1,s_2,\dots,s_b)$. Each vertical lozenge with displayed label $k$ is weighted by $f_{X,Y}(k)$. The other lozenges are weighed by $1$.}\label{Sshape2}
 \end{figure}
 
We first need the following two weighted enumeration formulas for dented semihexagons with removed unit triangles on the base.

Consider a semihexagon of side-lengths $a,b,a+b,b$ (in clockwise order, starting from the top). We remove $b$ up-pointing unit triangles along the base at the positions $1\le s_1<s_<\cdots<s_b\leq a+b$. Denote the region by $SH_{a,b}(s_1,s_2,\dots,s_b)$. We apply the weight assignment $\wt_{X,Y,z}$ to the region with the coordinate system $Oij$ whose $i$-axis runs along the base and whose $j$-axis passes through the midpoints of the two parallel sides as shown in Figure \ref{Sshape2}(a).

\begin{lem}[Theorem A.1 in \cite{LaiByunIntrusion2}]\label{lem1}
Let a and b be nonnegative integers, and let $\{s_1<\cdots<s_b\}\subseteq[a+b]$. Then
\begin{align*}
\M_{X,Y,z}(SH_{a,b}(s_1,s_2,\dots,s_b))=\prod_{1\le i <j\le b}f_{X,Y}(z + s_i + s_j - a - b - 1) \cdot \prod_{1\le i <j\le b}\frac{\br{s_j-s_i}}{\br{j-i}}.
\end{align*}
\end{lem}
 Lemma \ref{lem1} has been proved in  \cite[Theorem A.1]{LaiByunIntrusion2}.

For the second enumeration, reflect the region across a vertical axis and then rotate the reflected region counterclockwise by $120^\circ$.  Denote the resulting region by $SH'_{a,b}(s_1,s_2,\dots,s_b)$. We apply the weight assignment $\wt_{X,Y,z}$ to the region with the coordinate system shown in Figure  \ref{Sshape2}(b). In particular, the $i$-axis runs along the base of the region, and the $j$-axis passes through the its left endpoint.

\begin{lem}\label{lem2} Under the same hypotheses as in Lemma \ref{lem1}, we have
\begin{align*}
\M_{X,Y,z}(SH'_{a,b,z}(s_1,s_2,\dots,s_b))=\prod_{j=1}^{b} \prod_{r=1}^{s_j-j}f_{X,Y}(z+j+r-a-1) \cdot \prod_{1\le i <j\le b}\frac{\br{s_j-s_i}}{\br{j-i}}.
\end{align*}
\end{lem}

Lemma \ref{lem2} can be proved in the same way as Lemma \ref{lem1}. One could also prove these two lemmassing the Lindstr\"{o}m--Gessel--Viennot theorem \cite{Lind, GV, Stem3} and a Vandermonde-type determinant evaluation.

Before going to the proof, we need the following supplemental lemma.

\begin{lem}\label{extralem} Assume that $m,n,x,y$ are nonnegative integers. Set $m+n=N$, $A=\{a_1<a_2<\cdots<a_m\}$, $B=\{m+1,m+2,\dots,m+n\}$. Let 
\[A_x=\{N+1-a_m,\dots, N+1-a_1\}\cup \{N+1,\dots,N+x\}\]
\[A_y=\{N+1-a_m,\dots, N+1-a_1\}\cup \{N+1,\dots,N+y\}\]
 Then
\begin{equation}\label{exeq1}
\frac{E_{x,z}(B)}{E_{y,z}(B)} \cdot \frac{F_y(N)}{F_x(N)} \frac{D_x(A)D_x(B)}{D_y(A)D_y(B)}=\frac{\M_{Y,X,-z-2m}(SH'_{n,m+x}(A_x))}{\M_{Y,X,-z-2m}(SH'_{n,m+y}(A_y))}.
\end{equation}
\end{lem}
\begin{proof}
First, expand the quotient of E-factors
\[\frac{E_{x,z}(B)}{E_{y,z}(B)}=\frac{\prod_{k=1}^{x} \prod_{j=1}^{n}f_{X,Y}(z+m+j-k)}{\prod_{k=1}^{y} \prod_{j=1}^{n}f_{X,Y}(z+m+j-k)}\]
Since $D_x(B)=F_x(N)/F_x(m)$, the remaining B-factors simplify as follows
\begin{align}
\frac{F_y(N)}{D_y(B)}\frac{ D_x(B)}{F_x(N)}=\frac{F_y(N)}{\frac{F_y(N)}{F_y(m)}}\frac{ \frac{F_x(N)}{F_x(m)}}{F_x(N)}=\frac{F_y(m)}{F_x(m)}.
\end{align}
Thus, we write the expression on left-hand side of (\ref{exeq1}) as
\begin{align}\label{extralemeq1}
\frac{\prod_{i=1}^{x} \prod_{j=1}^{n}f_{X,Y}(z+m+j-i)}{\prod_{i=1}^{y} \prod_{j=1}^{n}f_{X,Y}(z+m+j-i)} \cdot \frac{D_x(A)F_y(m)}{D_y(A)F_x(m)}.
\end{align}

\bigskip

We now simplify the ratio of generating functions $\frac{\M_{Y,X,-2m-z}(SH'_{n,m+x}(A_x))}{\M_{Y,X,-2m-z}(SH'_{n,m+y}(A_y))}$.  Write
\begin{align*}
A_x&=\{a'_1<a'_2<\dots<a'_m<a'_{m+1}<\dots<a'_{m+x}\}\\
&= \{N+1-a_m,\dots, N+1-a_1\}\cup \{N+1,\dots,N+x\}=\overline{A} \cup B_x,\end{align*}
where \[\overline{A}=\{N+1-a_m,\dots, N+1-a_1\} \text{ and } B_x=\{N+1,\dots,N+x\}.\]
Thus, we have $a'_{i}=N+1-a_{m-i+1}$ for $1\leq i\leq m$, and $a'_{m+k}=N+k$,  for $k=1,2,\dots,x$.

We define the operation \[\Delta(S):=\prod_{1\le i<j \le l} \br{s_j-s_i},\] for any finite set of integers $S=\{s_1<s_2<\dots<s_l\}$. The cross version of $\Delta$ is defined by
\[\Delta(B,C):=\prod_{B\in C}\prod_{c\in C} \br{c-b},\] where every element of $B$ is smaller than every element of $C$. For such ordered disjoint sets,  we have
\[\Delta(B \cup C)=\Delta(B)\cdot \Delta(C) \cdot \Delta(B,C).\]

We have 
\[\M_{Y,X,-2m-z}(SH'_{n,m+x}(A_x))=\prod_{j=1}^{m+x}\prod^{a'_j-j}_{r=1}f_{Y,X}(-2m-z+j+r-n-1) \cdot \frac{\Delta_q(A_x)}{\Delta_q([m+x])}.\]

By the definition
\[f_{Y,X}(-2m-z+v)=f_{X,Y}(2m+z-v).\]
For every $k$ with $1\le k\le x$, we have $a'_{m+k}-(m+k)=(N+k)-(m+k)=n$. Thus, the contribution of term $a'_{m+k}$  to the weight product $\prod_{j=1}^{m+x}\prod^{a'_j-j}_{r=1}f_{Y,X}(-2m-z+j+r-n-1)$ is
\begin{align}
\prod_{r=1}^{n} f_{Y,X}(-2m-z+m+k+r-n-1)&=\prod_{r=1}^{n} f_{X,Y}(z+m+n+1-k-r)\notag\\
&=\prod_{j=1}^{n} f_{X,Y}(z+m+j-k).
\end{align}
The factors contributed by the first $m$ entries are independent of $x$, so they cancel in the ratio $\frac{\M_{Y,X,-2m-z}(SH'_{n,m+x}(A_x))}{\M_{Y,X,-2m-z}(SH'_{n,m+y}(A_y))}$. This implies that
\begin{align}\label{extralemeq2}\frac{\M_{Y,X,-2m-z}(SH'_{n,m+x}(A_x))}{\M_{Y,X,-2m-z}(SH'_{n,m+y}(A_y))}=\frac{\prod_{k=1}^{x} \prod_{j=1}^{n} f_{X,Y}(z+m+j-k) }{\prod_{k=1}^{y} \prod_{j=1}^{n} f_{X,Y}(z+m+j-k) }\cdot  \frac{\frac{\Delta(A_x)}{\Delta([m+x])}}{\frac{\Delta(A_y)}{\Delta([m+y])}}.
\end{align}

Set $C_x=\{m+1,\ldots,m+x\}$. Then $[m+x]=[m]\cup C_x$. Define $C_y$ similarly and have $[m+y]=[m]\cup C_y.$ We decompose the delta product of each of the sets $A_x,A_y,[m+x],[m+y]$ in terms of its components and have
\begin{align}
\displaystyle\frac{\frac{\Delta(A_x)}{\Delta([m+x])}}{\frac{\Delta(A_y)}{\Delta([m+y])}}=&\frac{\Delta(\overline{A}) \Delta(B_x)\Delta(\overline{A},B_x)}{\Delta(\overline{A}) \Delta(B_y)\Delta(\overline{A},B_y)}\cdot \frac{\Delta([m]) \Delta(C_y)\Delta([m],C_y)}{\Delta([m]) \Delta(C_x)\Delta([m],C_x)}\notag\\
&=\frac{\Delta(\overline{A},B_x)}{\Delta([m],C_x)}\frac{\Delta([m],C_y)}{\Delta(\overline{A},B_y)},
\end{align}
as $\Delta(B_x)=\Delta(C_x)$ and $\Delta(B_y)=\Delta(C_y)$.

We write 
\[\frac{\Delta(\overline{A},B_x)}{\Delta([m],C_x)}=\prod_{k=1}^{x} \prod_{i=1}^{m} \frac{\br{a'_{m+k}-a'_{m+1-i}}}{\br{(m+k)-(m+1-i)}}=\prod_{k=1}^{x}\prod_{i=1}^{m} \frac{\br{a_i+k-1}}{\br{i+k-1}}.\]
Similarly
\[\frac{\Delta(\overline{A},B_y)}{\Delta([m],C_y)}=\prod_{k=1}^{y}\prod_{i=1}^{m} \frac{\br{a_i+k-1}}{\br{i+k-1}}.\]
Dividing the two identities, we have
\begin{align}\label{extralemeq3}
\frac{\Delta(\overline{A},B_x)}{\Delta([m],C_x)}\frac{\Delta([m],C_y)}{\Delta(\overline{A},B_y)}&=\frac{\prod_{k=1}^{x}\prod_{i=1}^{m} \frac{\br{a_i+k-1}}{\br{i+k-1}}}{\prod_{k=1}^{y}\prod_{i=1}^{m} \frac{\br{a_i+k-1}}{\br{i+k-1}}}\notag\\
&=\frac{\prod_{i=1}^{m}\frac{\br{a_i}_{x}}{\br{i}_x}}{\prod_{i=1}^{m}\frac{\br{a_i}_{y}}{\br{i}_y}}=\frac{D_x(A)F_y(m)}{D_y(A)F_x(m)}.
 \end{align}
By (\ref{extralemeq1}), (\ref{extralemeq2}), and (\ref{extralemeq3}), both sides of (\ref{extralemeq1}) equal
\[\frac{\prod_{k=1}^{x} \prod_{j=1}^{n}f_{X,Y}(z+m+j-k)}{\prod_{k=1}^{y} \prod_{j=1}^{n}f_{X,Y}(z+m+j-k)} \cdot \prod_{i=1}^m\frac{D_x(A)F_y(m)}{D_y(A)F_x(m)}.\]
\end{proof}

\begin{figure}[ht]\centering
          \includegraphics[width=14cm]{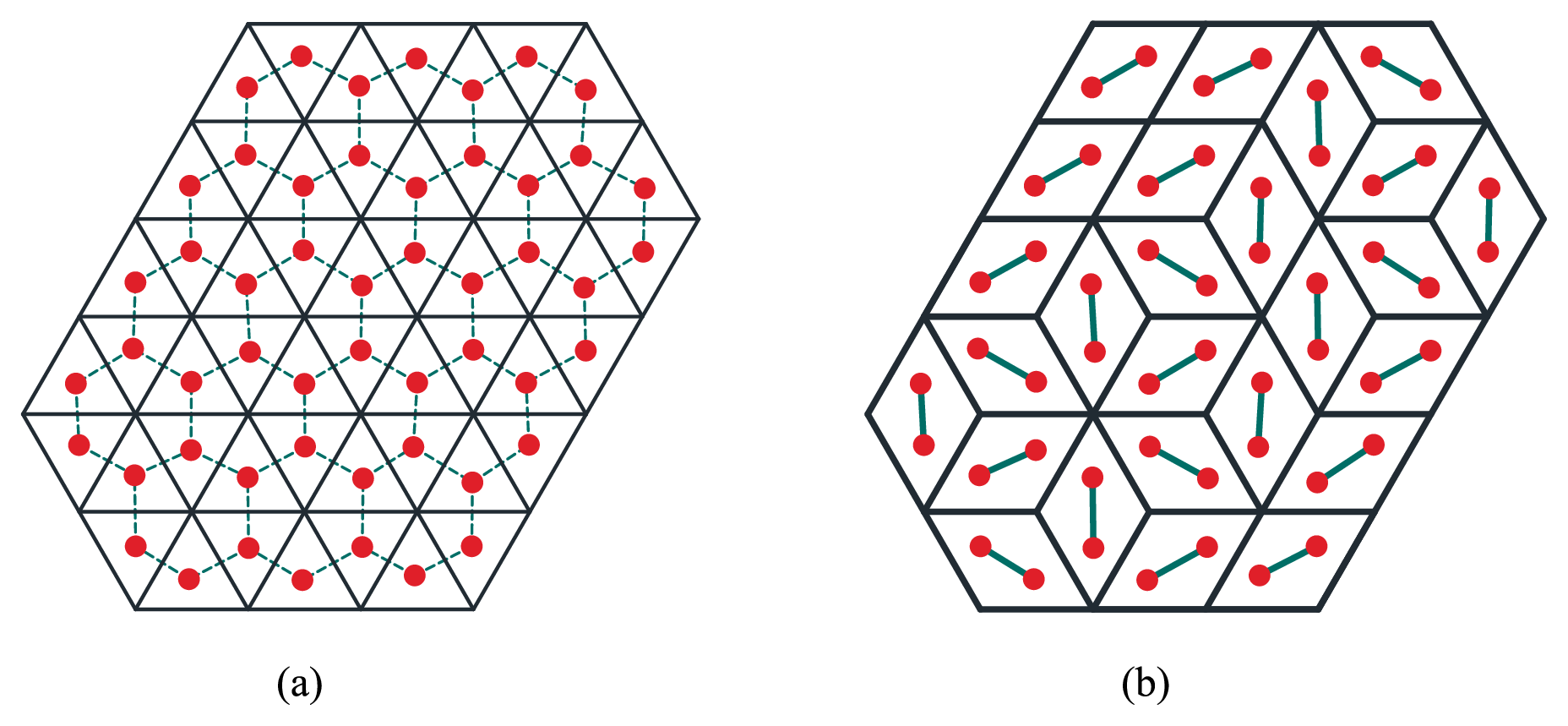}
 \caption{(a) The hexagon with side lengths $4,$ $3,$ $2,$ $4,$ $3,$ $2$ and its dual graph. Red dots represent
vertices, and dashed segments represent edges. (b) A lozenge tiling and the corresponding perfect
matching}\label{dualgraph}
 \end{figure}

 The dual graph $G$ of a region $R$ on the triangular lattice is the graph whose vertices are the unit triangles in $R$ and whose edges connect exactly the two unit triangles sharing an edge. Each edge in $G$ is weighted by the same weight as the corresponding lozenge in $R$. A \textit{perfect matching} in a graph  $G$ is a collection of vertex-disjoint edges that covers all vertices of the graph. The \textit{weight} of a perfect matching is the product of weights of its edges. This way, we get a weight-preserve bijection between perfect matchings of $G$ and tilings of $R$. See Figure \ref{dualgraph} for an example. We use the notation $\M(G)$ for the sum of weights of its perfect matchings.

\begin{lem}[Theorem 5.4 in \cite{Kuo}]\label{lemkuo}
Let $G=(V_1,V_2,E)$ be a weighted plane bipartite graph in which $|V_1|=|V_2|+2$. Let vertices $u,v,w,s$ appear on a face of $G$, in that order. If $u,v,w,s\in V_1$, then 
\begin{align}\label{kuoeq}
\M(G-\{u,w\})\M(G-\{v,s\})=&\M(G-\{u,v\})\M(G-\{w,s\})\notag\\
&\quad \quad +\M(G-\{u,s\})\M(G-\{v,w\}),
\end{align}
where $G-\{u,w\}$ denotes the graph obtained from $G$ by removing the two vertices $u$ and $w$ together with all edges incident to $u$ or $w$, and where other graphs are defined similarly.
\end{lem}

\bigskip

We are now ready to move the the proof of Lemma \ref{shufflingthm2}.
\begin{proof}[Proof of Lemma \ref{shufflingthm2}]
Set $N=m+n$. Assume that $A$ and $B$ are two subsets of $[N]$ with $A=\{a_1<a_2<\cdots<a_m\}$ and $|B|=\{b_1<b_2<\cdots<b_n\}$.
We only consider the case when our regions $S'_{N, x}(A,B)$ and $S'_{N,y}(A,B)$ have a tiling. By Lemma \ref{lem0}, it is equivalent that 
\[|A\cap[k]|+|B\cap [k]|\leq k,\]
for every $k$ with $1\le k \le m+n$.
It suffices to treat $x\ge y$; the other case follows by interchanging $x$ and $y$ and inverting the ratio. 

We rewrite the identity as
\begin{align}\label{proofeq1}
\M(S'_{N,x}(A,B))=g_{x,y,z}(A,B)\cdot \M(S'_{N,y}(A,B)),
\end{align}
where 
\begin{align}\label{proofeq2x}
g_{x,y,z}(A,B)=\frac{E_{x,z}(B)}{E_{y,z}(B)} \cdot \frac{F_y(|A|+|B|)}{F_x(|A|+|B|)} \frac{D_x(A)D_x(B)}{D_y(A)D_y(B)}.
\end{align}
To clarify individual factors, one may introduce the temporary abbreviations
\[
Q_{x,y}(s)=\frac{\langle s\rangle_x}{\langle s\rangle_y},\qquad P_{x,y,z}(s)=\frac{\prod_{i=1}^{x}f_{X,Y}(z+s-i)}{\prod_{i=1}^{y}f_{X,Y}(z+s-i)}.
\]
Then the same formula reads
\begin{align}\label{alternativeg}
g_{x,y,z}(A,B)=\frac{F_y(|A|+|B|)}{F_x(|A|+|B|)}\prod_{a\in A}Q_{x,y}(a)\prod_{b\in B}Q_{x,y}(b)P_{x,y,z}(b).
\end{align}
For $s>1$, we have the following useful identities for $Q$- and $P$-products
\[
Q_{x+1,y+1}(s-1)=Q_{x,y}(s),
\]
\[
P_{x+1,y+1,z+2}(s-1)=P_{x,y,z}(s).
\]

We prove (\ref{proofeq1}) by induction on the statistic $p:=m+n+t=N+t$, where
 \[t:=\sum_{i=1}^{n} ((m+i)-b_i)=mn+\frac{n(n+1)}{2}-\sum_{b\in B} b=mn+\frac{n(n+1)}{2}-\sum_{i=1}^{n} b_i.\] 
 The statistic t measures the total displacement of the $B$-dents from the terminal positions $m+1,\ldots ,m+n$. As $B$ is an $n$-subset of $[m+n]$, we have $b_i\leq m+i$, thus $t\geq 0$.  When $t=0$, all these removed unit triangles are clustered at the lower-right corner, i.e. $B=\{m+1,m+2,\dots,N\}$ (see Figure \ref{basecase2}(a)). The base cases are the situations when at least one of the parameters $m,n,t$ is equal to $0$.

   \begin{figure}[ht]\centering
         \includegraphics[width=15cm]{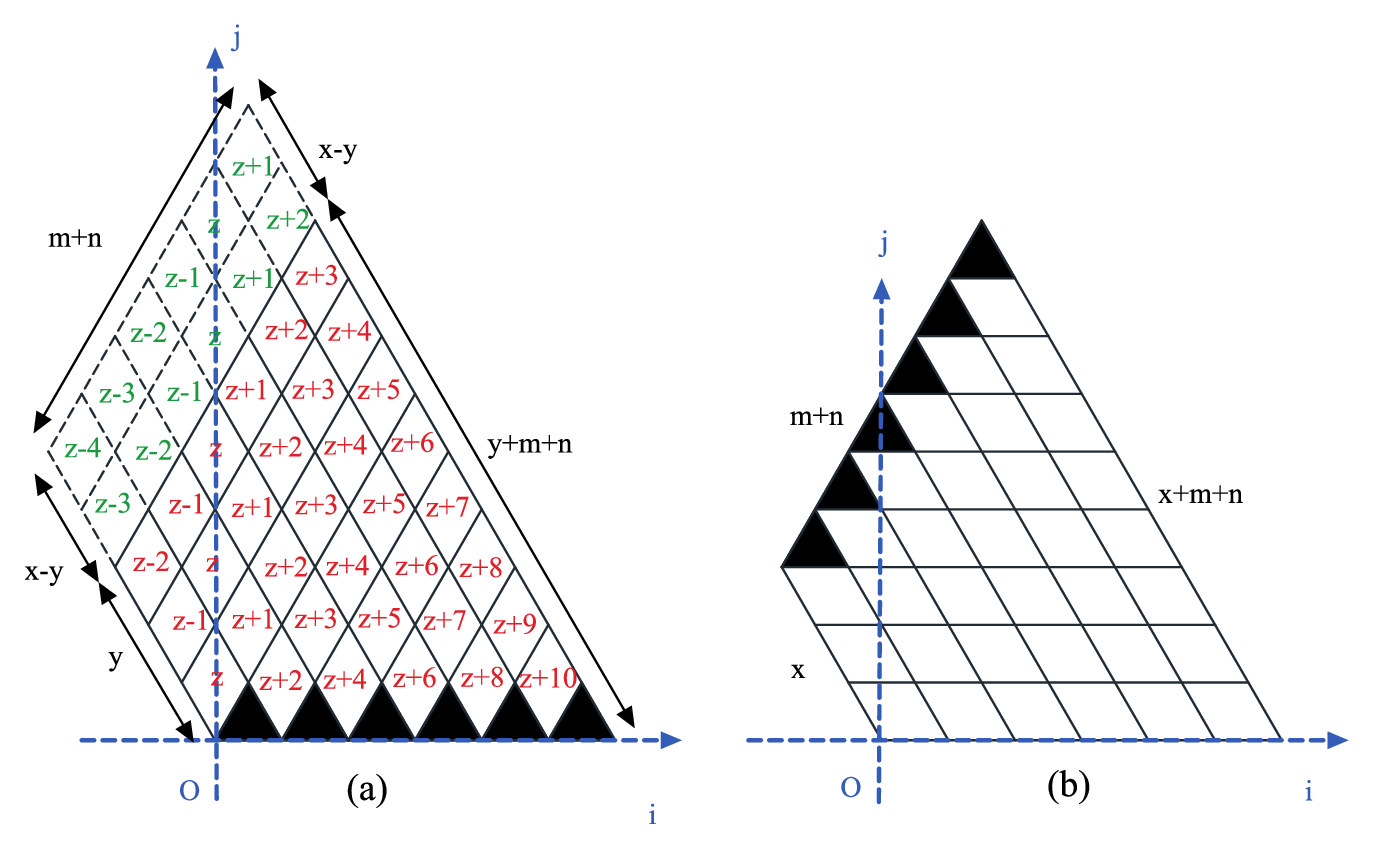}
 \caption{The base cases of the induction. (a) The unique tilings of $S'_{n,x}(\varnothing,[n])$ and $S'_{n,y}(\varnothing,[n])$ superimposed for $x\ge y$; the additional lozenges have dashed outlines and green labels. (b) The unique tiling of $S'_{m,x}([m],\varnothing)$, whose lozenges all have weight 1. The labels in (a) are the arguments of $f_{X,Y}$.} \label{basecase}
 \end{figure}

If $m=0$, then $A=\emptyset$ and $B=[n]$ and each of the semi-hexagons $S'_{N,x}(A,B)$ and $S'_{N,y}(A,B)$ has exactly one tiling consisting entirely of the vertical lozenges. We superimpose two regions in Figure \ref{basecase}(a). This implies that
$\frac{\M_{X,Y,z}(S'_{N,x}(A,B))}{\M_{X,Y,z}(S'_{N,y}(A,B))}$
will be the product of weights of the lozenges with dashed boundary with green labels. Counting these lozenges explicitly, we have
\begin{align}
\frac{\M_{X,Y,z}(S'_{N,x}(A,B))}{\M_{X,Y,z}(S'_{N,y}(A,B))} &=\prod_{i=1}^{x-y}\prod_{j=1}^{l} f_{X,Y}(z-y-i+j)\notag\\
&=\prod_{r=y+1}^{x}\prod_{j=1}^{l} f_{X,Y}(z+j-r)=\frac{\prod_{r=1}^{x}\prod_{j=1}^{l} f_{X,Y}(z+j-r)}{\prod_{r=1}^{y}\prod_{j=1}^{l} f_{X,Y}(z+j-r)}.
\end{align}
The final expression is exactly $g_{x,y,z}(\emptyset, [n])$ in (\ref{proofeq2x}). The symmetric $q$-integer factors cancel, leaving only the $f_{X,Y}$ factors.. Then, the identity (\ref{proofeq1}) holds when $m=0$.

If $n=0$, then $A=[m]$ and $B=\emptyset$ and each of the semi-hexagons $S'_{N,x}(A,B)$ and $S'_{N,y}(A,B)$ has exactly one tiling consisting of all the left-tiling lozenges, all of wight $1$. See the picture Figure \ref{basecase}(b) for the region $S'_{N,x}([m],\emptyset)$. It means that $\frac{\M_{X,Y,z}(S'_{N,x}(A,B))}{\M_{X,Y,z}(S'_{N,y}(A,B))}=1$. In this case, 
\[g_{x,y,z}([m],\emptyset)=\frac{F_y(m)}{F_x(m)}\frac{D_x([m])D_x(\emptyset)}{D_y([m])D_y(\emptyset)}=\frac{F_y(m)}{F_x(m)}\frac{F_x(m)}{F_y(m)}=1.\]
 Thus, the identity (\ref{proofeq1}) holds when $n=0$.

If $t=0$, then \[\M_{X,Y,z}(S'_{N,x}(A,B))=V\cdot \M_{X,Y,z}(P),\] where $P$ is the pentagonal region obtained from $S'_{N,x}(A,B)$ by removing the vertical lozenges in the triangle in the lower-right corner and $V$ is the product of  weights of these forced lozenges (see Figure \ref{basecase2}(a)). The region $P$ can also be obtained from the trapezoidal region $H$ in Figure \ref{basecase2}(a) by removing forced lozenges. These lozenges all have weight $1$. Thus 
\[\M_{X,Y,z}(P)=\M_{X,Y,z}(H).\]
The region $H$ is obtained by translating horizontally the semihexagon $SH'_{n,m+x}(A_x)$ in Lemma \ref{lem2} by $-2m$ units (i.e., by $2m$ units to the negative $i$-direction) and reflecting the resulting region over the $j$-axis. Denote by $S^{*}$ the intermediate region. The reflecting property (P2) gives
\[\M_{X,Y,z}(H)=\M_{Y,X,-z}(S^{*}),\]
and  the translating property (P3) gives
\[\M_{Y,X,-z}(S^{*})=\M_{Y,X,-2m-z}(SH'_{n,m+x}(A_x)).\]
Thus, by the four identities above, we have
\[\M_{X,Y,z}(S'_{m+x,x}(A,B))=V\cdot\M_{Y,X,-2m-z}(SH'_{n,m+x}(A_x)),\]  where $F$ is the product of weight of forced lozenges. Similarly, we have
\[\M_{X,Y,z}(S'_{m+y,y}(A,B))=V\cdot\M_{Y,X,-2m-z}(SH'_{n,m_y}(A_y)).\] 
Thus, we have
\[\frac{\M_{X,Y,z}(S'_{N,x}(A,B))}{\M_{X,Y,z}(S'_{N,y}(A,B))}=\frac{M_{Y,X,-2m-z}(SH'_{n,m+x}(A_x))}{\M_{Y,X,-2m-z}(SH'_{n,m+y}(A_y))}\]
and  the identity (\ref{proofeq1}) follows from Lemma \ref{lem2}.

\begin{figure}[ht]\centering
         \includegraphics[width=15cm]{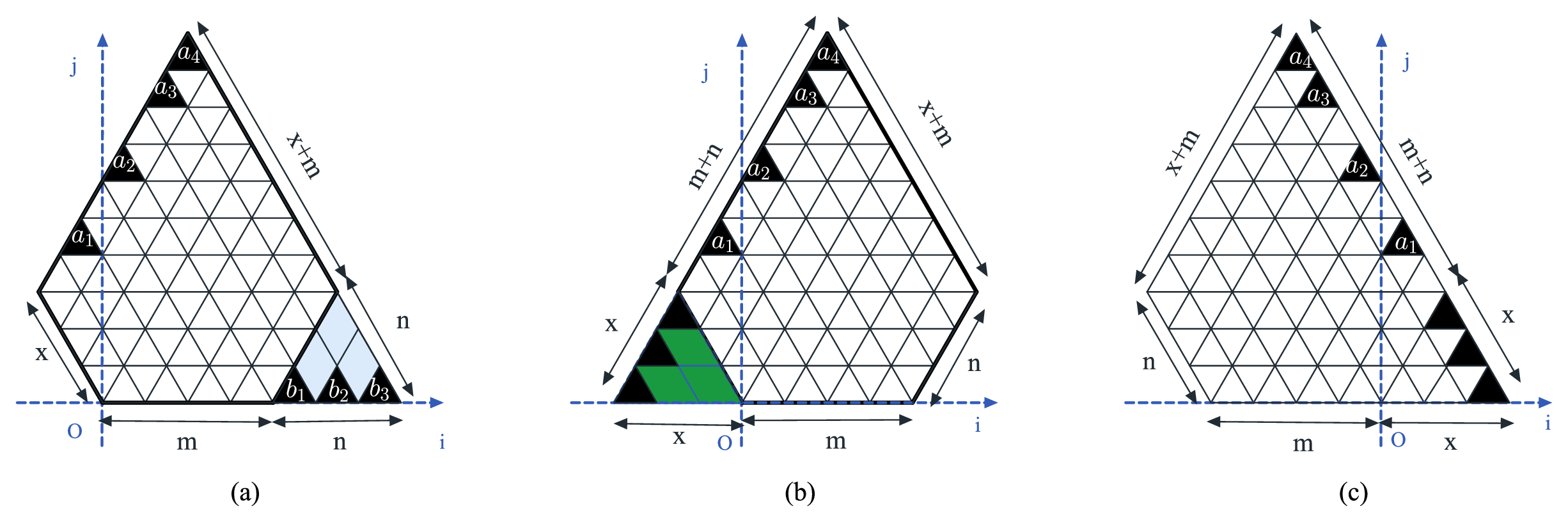}
 \caption{(a) Removing forced vertical lozenges from $S'_{N,x}(A,B)$ when $B=\{m+1,\ldots ,N\}$ leaves the pentagonal region P. (b) Removing forced lozenges of weight 1 from H gives the same region P. (c) Translating $S'_{n,m+x}(A_x)$ by $-2m$ horizontal coordinate units and then reflecting across the $j$-axis gives $H$.} \label{basecase2}
 \end{figure}
 
\begin{figure}[ht]\centering
         \includegraphics[width=12cm]{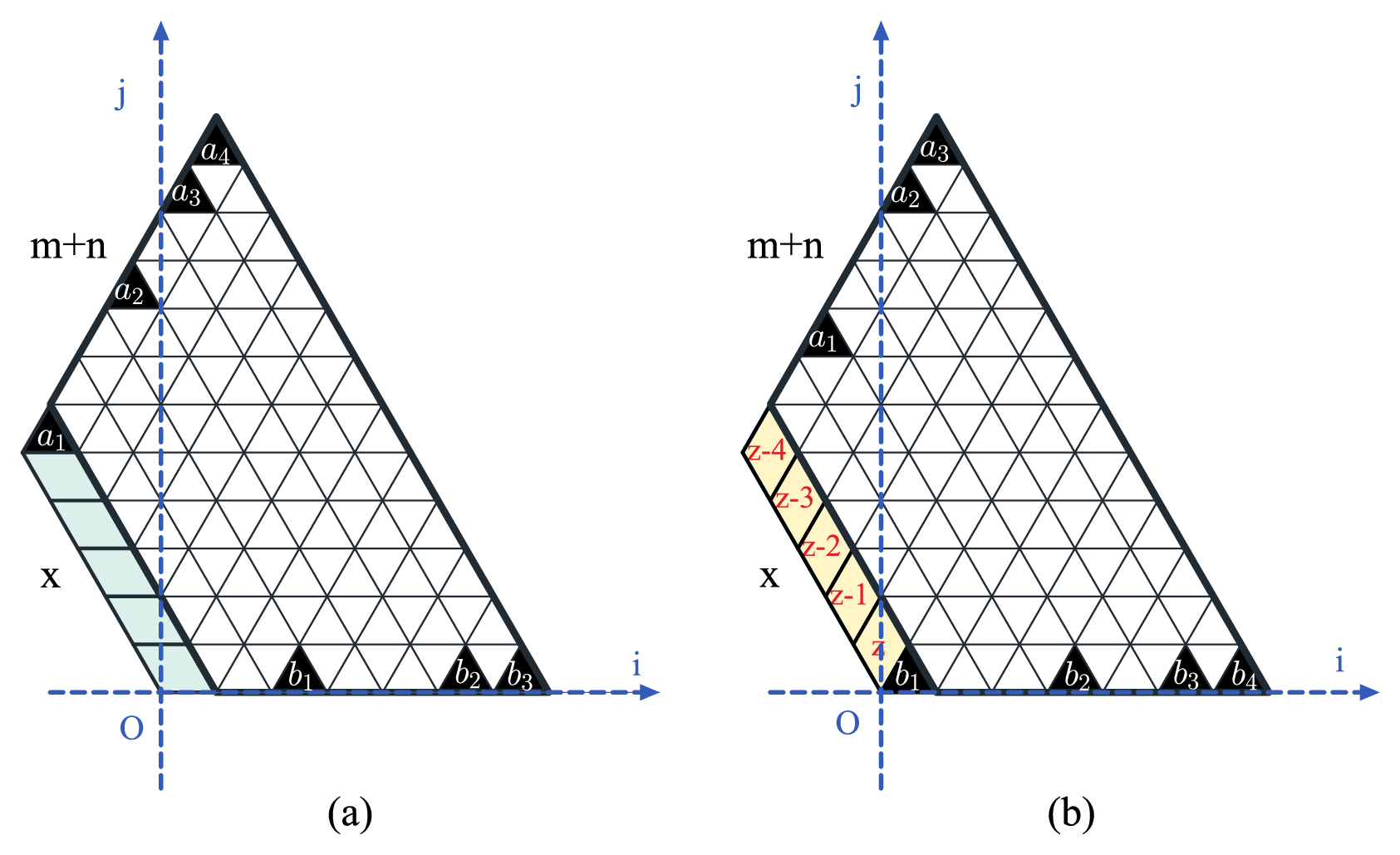}
 \caption{Picture (a): the case $a_1=1$; the removal of forced lozenges (shaded ones) from the region $S'_{N,x}(A,B)$ gives a region with smaller $p$-statistic. Picture (b): the case $a_1=2$; the removal of forced lozenges (shaded ones) from the region $S'_{N,x}(A,B)$ gives a region with smaller $p$-statistic.} \label{basecase3}
 \end{figure}

For the induction step, we assume that $m,n,t$ are all positive and that (\ref{proofeq1}) holds for any pair of semihexagons whose $p$-statistic is strictly less than $m+n+t$.

For a finite set of integers $A=\{a_1,a_2,\dots,a_k\}$ we use the shorthand notation \[T(A)=A-1=\{a_1-1,a_2-1,\dots,a_k-1\}\] for the translate of $A$ by $-1$ in the rest of the proof.
By Lemma \ref{lem0}, at least one of $a_1$ and $b_1$ is strictly greater than $1$. If $a_1=1$, then $b_1>1$ and one can remove forced left-tilting lozenges along the southwest side of region $S'_{N,x}(A,B)$ to obtain a region $\widehat{S}$ with smaller $p$-statistic (see Figure \ref{basecase3}(a)). The region $\widehat{S}$ is obtained from the region $S'_{N-1,x+1}(T(A\setminus\{1\}),T(B))$ by shifting $2$ units to the right.  We recall that one horizontal unit is half the side length of a unit triangle. These new regions have the same tiling generating functions as the original as the removed forced lozenges all have weight $1$. We then have
\[\M_{X,Y,z}(S'_{N,x}(A,B))=\M_{X,Y,z}(\widehat{S})= \M_{X,Y,z+2}(S'_{N-1,x+1}(T(A\setminus\{1\}),T(B))).\]
It is easy to verify that the $n$-statistic of the new region is $n$. Its $m$-statistic is $m-1$, and its $t$-statistic is 
\[\sum_{i=1}^{m} ((m-1)+ i- (b_i-1))=\sum_{i=1}^{m} (m+ i- b_i)=t.\] Thus the $p$-statistic of the new region is $p-1<p$.
Similarly, we have 
\[\M_{X,Y,z}(S'_{N,y}(A,B))= \M_{X,Y,z+2}(S'_{N-1,y+1}(T(A\setminus\{1\}),T(B))).\]
Since the $m,n,p$-statistic are independent of $x$, the region on the right-hand side has its $p$-statistic $p-1$, which is less then the $p$-statistic of the region on the left-hand side. Dividing the two identities and applying the induction hypothesis gives
\begin{align}
\frac{\M_{X,Y,z}(S'_{N,x}(A,B))}{\M_{X,Y,z}(S'_{N,y}(A,B))}&=\frac{M_{X,Y,z+2}(S'_{N-1,x+1}(T(A\setminus\{1\}),T(B)))}{M_{X,Y,z+2}(S'_{N-1,y+1}(T(A\setminus\{1\}),T(B)))}\notag\\
&=g_{x+1,y+1,z+2}(T(A\setminus\{1\}),T(B)).
\end{align}
To verify (\ref{proofeq1}), we need to show that
\[g_{x,y,z}(A,B)=g_{x+1,y+1,z+2}(T(A\setminus\{1\}),T(B)),\]
provided that $a_1=1$ and $b_1>1$.

To clarify individual factors, one may introduce the temporary abbreviations
\[
Q_{x,y}(s)=\frac{\langle s\rangle_x}{\langle s\rangle_y},\qquad P_{x,y,z}(s)=\frac{\prod_{i=1}^{x}f_{X,Y}(z+s-i)}{\prod_{i=1}^{y}f_{X,Y}(z+s-i)}.
\]
We can now rewrite function $g_{x,y,z}$ as
\begin{align}\label{alternativeg}
g_{x,y,z}(A,B)=\frac{F_y(|A|+|B|)}{F_x(|A|+|B|)}\prod_{a\in A}Q_{x,y}(a)\prod_{b\in B}Q_{x,y}(b)P_{x,y,z}(b).
\end{align}
For $s>1$, we have the following useful identities for $Q$- and $P$-products
\begin{align}\label{reduceeq1}
Q_{x+1,y+1}(s-1)=Q_{x,y}(s),
\end{align}
\begin{align}\label{reduceeq2}
P_{x+1,y+1,z+2}(s-1)=P_{x,y,z}(s).
\end{align}
We also have for any positive integer $k$
\begin{align}\label{reduceeq3}
\frac{F_{y+1}(k-1)}{F_{x+1}(k-1)}=\frac{F_y(k)}{F_x(k)}Q_{x,y}(1).
\end{align}

We simplify the ratio of the two $g$-functions using formula (\ref{alternativeg}) and identities (\ref{reduceeq1}), (\ref{reduceeq2}), and (\ref{reduceeq3})
\begin{align}
&\frac{g_{x+1,y+1,z+2}(T(A\setminus\{1\}),T(B))}{g_{x,y,z}(A,B)}\notag\\
&=\frac{F_x(N)F_{y+1}(N-1)}{F_{x+1}(N-1)F_{y}(N)} \cdot \frac{\prod_{a\in T(A\setminus\{1\})} Q_{x+1,y+1}(a)}{\prod_{a\in A} Q_{x,y}(a)} \cdot \frac{\prod_{b\in T(B)}Q_{x+1,y+1}(b) P_{x+1,y+1,z+2}(b)}{\prod_{b\in B}Q_{x,y}(b)P_{x,y,z}(b)}\notag\\
&= Q_{x,y}(1) \cdot  \frac{\prod_{a\in A\setminus\{1\}} Q_{x+1,y+1}(a-1)}{\prod_{a\in A} Q_{x,y}(a)} \cdot \frac{\prod_{b \in B}Q_{x+1,y+1}(b-1) P_{x+1,y+1,z+2}(b-1)}{\prod_{b\in B}Q_{x,y}(b)P_{x,y,z}(b)}\notag\\
&= Q_{x,y}(1) \cdot \frac{1}{Q_{x,y}(1)} \cdot 1=1.
\end{align}
Similarly, for the case $b_1=1$ and $a_1>1$, removal of forced lozenges from $S'_{N,x}(A,B)$ as in (see Figure \ref{basecase3}(c)) gives a shifted version of $S'_{N-1,x+1}(T(A),T(B\setminus\{1\}))$. We have
\[\M_{X,Y,z}(S'_{N,x}(A,B))= \prod_{i=1}^{x}f_{X,Y}(z+1-i) \cdot \M_{X,Y,z+2}(S'_{N-1,x+1}(T(A),T(B\setminus\{1\}))),\]
where $\prod_{i=1}^{x}f_{X,Y}(i)$ is the product of weights of forced lozenges. 
 Similarly, we have
\[\M_{X,Y,z}(S'_{N,y}(A,B))= \prod_{i=1}^{y}f_{X,Y}(z+1-i) \cdot \M_{X,Y,z+2}(S'_{N-1,y+1}(T(A),T(B\setminus\{1\}))).\]
We can verify that  the regions on the right-hand sides of the above identity have $m$-statistic $m$, $n$-statistic $n-1$, and the $t$-statistic given by
\[\sum_{i=1}^{n-1} (m+i-(b_{i+1}-1))=\sum_{j=2}^{n} (m+j-b_j)=t-(m+1-b_1)=t-m<t.\]
Thus they have $p$-statistic less than that of the corresponding regions on the left-hand sides.
Dividing the two identities and applying the induction hypothesis gives
\begin{align}
\frac{\M_{X,Y,z}(S'_{N,x}(A,B))}{\M_{X,Y,z}(S'_{N,y}(A,B))}&=\frac{\prod_{i=1}^{x}f_{X,Y}(z+i-1)}{\prod_{i=1}^{y}f_{X,Y}(z+i-1)}\frac{M_{X,Y,z+2}(S'_{N-1,x+1}(T(A),T(B\setminus\{1\})))}{M_{X,Y,z+2}(S'_{N-1,y+1}(T(A),T(B\setminus\{1\})))}\notag\\
&=\frac{\prod_{i=1}^{x}f_{X,Y}(z+1-i)}{\prod_{i=1}^{y}f_{X,Y}(z+1-i)}g_{x+1,y+1,z+2}(T(A),T(B\setminus\{1\}))).
\end{align}
To verify (\ref{proofeq1}) in this case, we need to show that
\[g_{x,y,z}(A,B)=\frac{\prod_{i=1}^{x}f_{X,Y}(z+1-i)}{\prod_{i=1}^{y}f_{X,Y}(z+1-i)} g_{x+1,y+1,z+2}(T(A),T(B\setminus\{1\})).\]
provided that $a_1>1$ and $b_1=1$.
Indeed, we use formula (\ref{alternativeg}) and identities (\ref{reduceeq1}), (\ref{reduceeq2}), and (\ref{reduceeq3}) to simplify the ratio of the two $g$-functions above
\begin{align}
&\frac{g_{x+1,y+1,z+2}(T(A),T(B\setminus\{1\}))}{g_{x,y,z}(A,B)}\notag\\
&=\frac{F_x(N)F_{y+1}(N-1)}{F_{x+1}(N-1)F_y(N)} \cdot \frac{\prod_{a\in T(A)} Q_{x+1,y+1}(a)}{\prod_{a\in A} Q_{x,y}(a)} \cdot \frac{\prod_{b\in T(B\setminus\{1\})}Q_{x+1,y+1}(b) P_{x+1,y+1,z+2}(b)}{\prod_{b\in B}Q_{x,y}(b)P_{x,y,z}(b)} \notag\\
&=Q_{x,y}(1) \cdot  \frac{\prod_{a\in A} Q_{x+1,y+1}(a-1)}{\prod_{a\in A} Q_{x,y}(a)} \cdot \frac{\prod_{b \in B\setminus \{1\}}Q_{x+1,y+1}(b-1) P_{x+1,y+1,z+2}(b-1)}{\prod_{b\in B}Q_{x,y}(b)P_{x,y,z}(b)}\notag\\
&= Q_{x,y}(1)\cdot 1\cdot \frac{1}{Q_{x,y}(1)P_{x,y,z}(1)}=\frac{1}{P_{x,y,z}(1)}=\frac{\prod_{i=1}^{y} f_{X,Y}(z+1-i)}{\prod_{i=1}^{x} f_{X,Y}(z+1-i)}.
\end{align}
Then  (\ref{proofeq1}) follows.

Thus, we can assume that  $a_1,b_1>1$. We will show that the expressions on both sides of (\ref{proofeq1}) satisfy the same recurrence.

     \begin{figure}[ht]\centering
         \includegraphics[width=8cm]{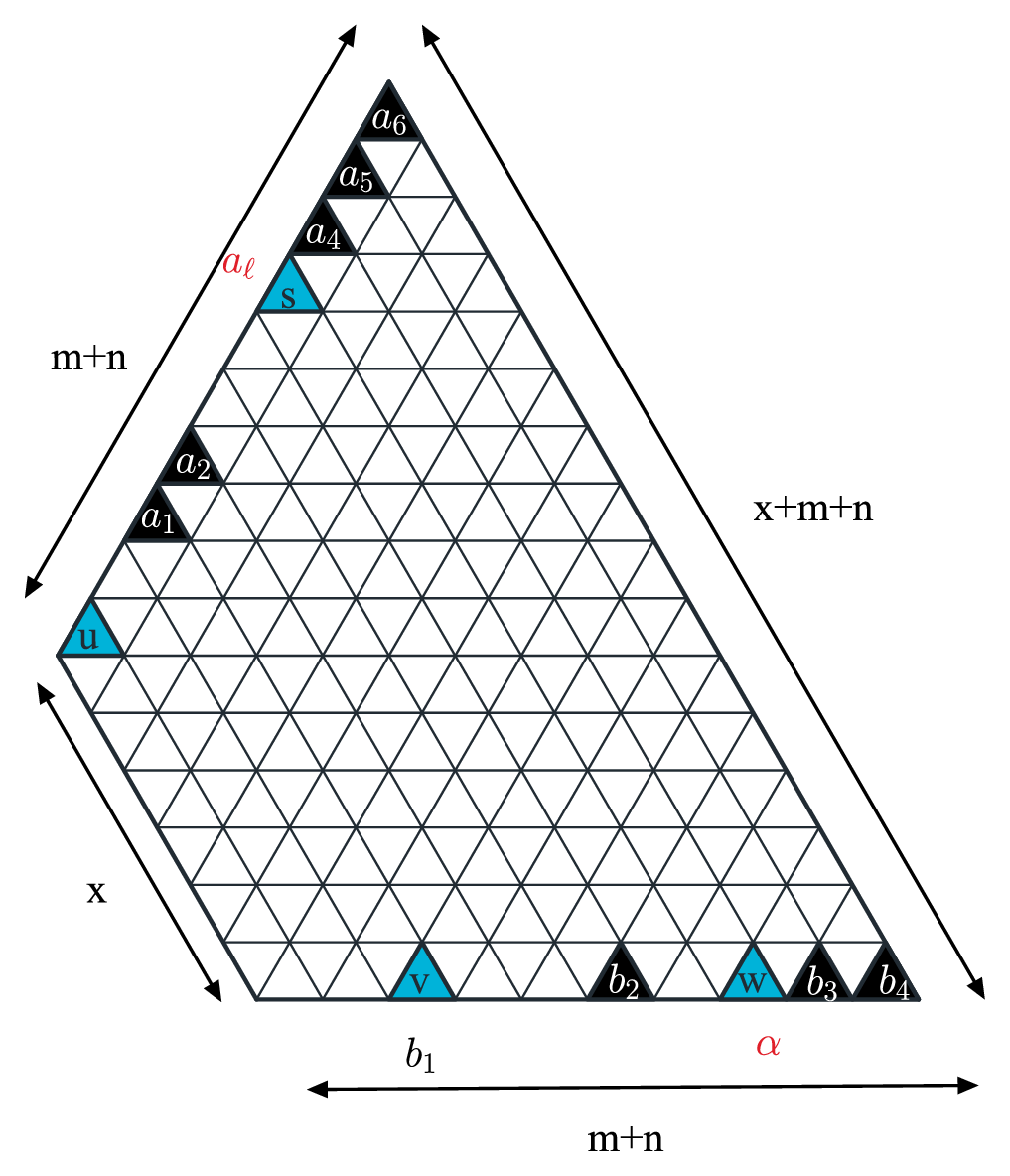}
 \caption{How we apply Kuo's condensation Lemma \ref{lemkuo} to the dual graph of the region $R^+$ obtained from $S_{N,x}(A,B)$ by restoring the removed unit triangles  at the positions $a_{\ell}$ and $b_1$. The vertices $u,v,w,s$ correspond with the shaded unit triangles of the same label. In this picture $\ell=3$, $b_1=3$ and $\alpha=8$.} \label{kuo1}
 \end{figure}
 
\begin{figure}\centering
         \includegraphics[width=12cm]{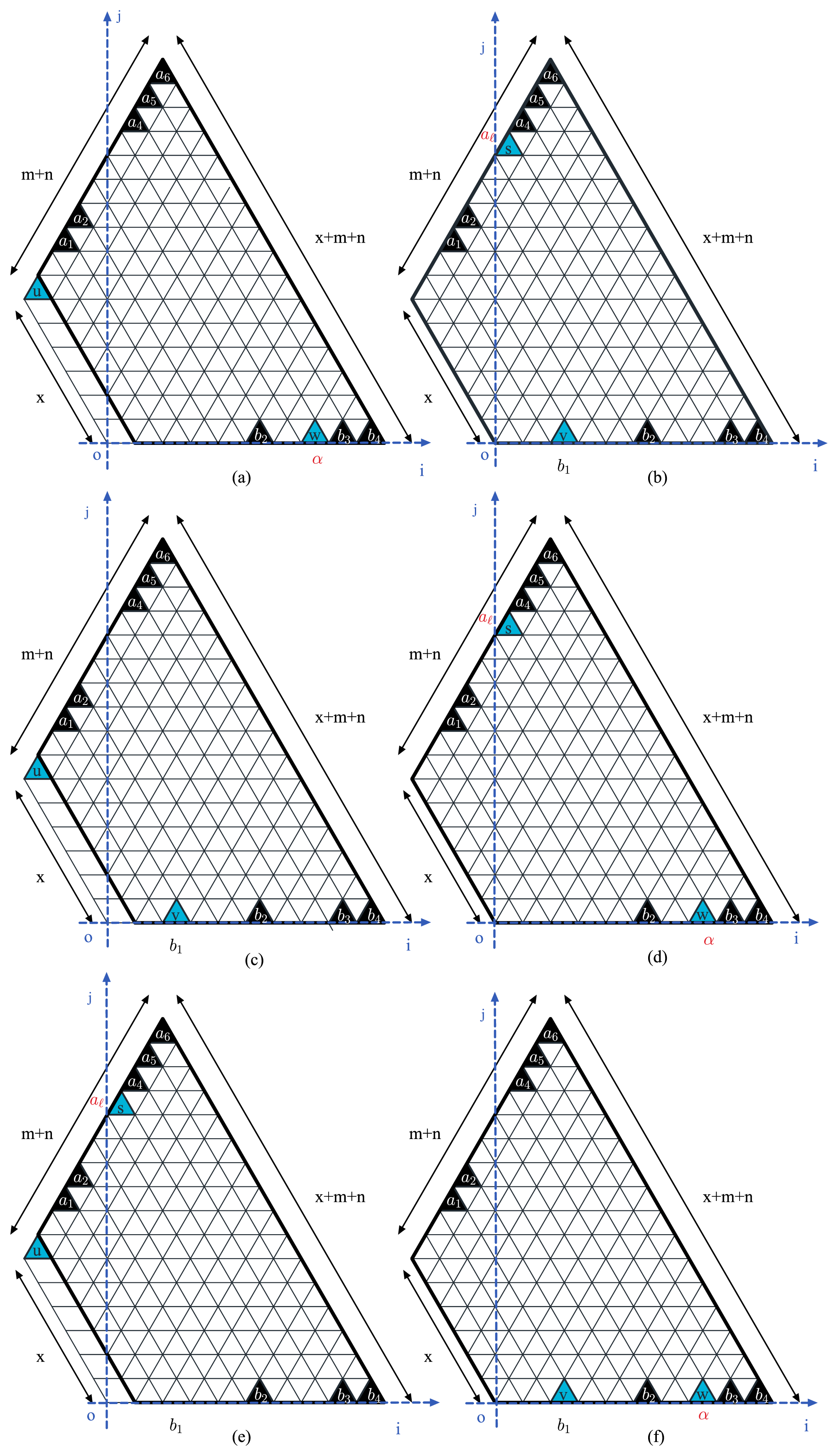}
 \caption{The six regions corresponding to $G-\{u,w\}$, $G-\{v,s\}$, $G-\{u,v\}$, $G-\{w,s\}$, $G-\{u,s\}$, and $G-\{v,w\}$, respectively, in pictures (a)--(f). Shaded triangles mark the deleted vertices. The bold contours enclose the regions remaining after the indicated forced lozenges are removed.} \label{kuo2}
 \end{figure}
 
 We obtain a recurrence for  the expression on left-hand side of (\ref{proofeq1})  by using Kuo's condensation in Lemma \ref{lemkuo}. Assume that $\ell$ is the largest index such that $a_{\ell}-1 \notin A$. We consider the region $R$ obtained from $S_{N,x}(A,B)$ by restoring the removed unit triangles  at the positions $a_{\ell}$ and $b_1$ (see Figure \ref{kuo1} for an example).  We now apply Kuo's condensation for the dual graph $G$ of $R^+$ with the four vertices $u,v,w,s$ are at the positions of the shaded up-pointing triangles with the same labels in  Figure \ref{kuo1}. Since the four triangles $u,v,w,s$ appear in counterclockwise order on the boundary of $R^+$, the vertices $u,v,w,s$ appear in a cyclic order are on the outer face of the dual graph $G$. The restoring two up-pointing dents in $R^+$ (at the positions $a_{\ell}$ and $b_1$) increases the corresponding vertex class in $G$ by $2$.
 
 In particular, the triangles with labels $v$ and $s$ are the positions $b_1$ and $a_3$, respectively. The $u$-triangle is the up-pointing unit triangle in the upper-left corner, and the $w$-triangle is the last non-remove triangle on the base of the region. Let $\alpha$ be the position of the $w$-triangle, so $\alpha=\max([N]\setminus B)$. 
 
  We note that our configuration as in Figure \label{kuo1} is needed for the induction. In particular, we must have the position $\alpha$ of $w$ stays strictly to the right of the position $b_1$ of $v$. Indeed, we always have $b_1\leq m< \alpha$ when $t>0$. Indeed, if $b_1\geq m+1$, then the $n$-element set $B$ must be the terminal block, so $t=0$. If $\alpha\leq m$, every position $m+1,\dots,N$ must be in $B$, again $t=0$.

We first consider the region corresponding with the graph $G-\{u,w\}$. It is illustrated in Figure \ref{kuo2}(a). The removal of $u$-triangle yields forced left-tilting lozenges running along the southwest side. This lozenge removal preserves the tiling generating function as the weights of the lozenges are all $1$.  The resulting region (enclosed by the bold contour) has the same tiling generating function as the original region. The resulting region is congruent to $\mathcal{R}_{u,w}:=S_{N-1,x+1}(T(A\setminus \{a_{\ell}\}); T((B\cup\{\alpha\})\setminus\{b_1\}))$. The coordinate systems do not yet agree. Translate the region $2$ units to the left so that the $j$-axis passes through the lower-left corner.  Thus we have
\[\M_{X,Y,z}(G-\{u,w\})=\M_{X,Y,z+2}(S'_{N-1,x+1}(T(A\setminus \{a_{\ell}\}); T((B\cup\{\alpha\})\setminus\{b_1\}))).\]
The region with respect to $G-\{v,s\}$ is exactly the region $\mathcal{R}=S'_{x}(A,B)$. Thus, we have
\[\M_{X,Y,z}(G-\{v,s\})=\M_{X,Y,z}(S'_{N,x}(A, B)).\]
In this case we do not need to perform any further translation to correct the weight assignment (see Figure \ref{kuo2}(b)).

Similarly, as shown in Figures \ref{kuo2}(c)--(f),
\[\M_{X,Y,z}(G-\{u,v\})=\M_{X,Y,z+2}(S'_{N-1,x+1}(T(A\setminus \{a_{\ell}\}); T(B))),\]
\[\M_{X,Y,z}(G-\{w,s\})=\M_{X,Y,z}(S'_{N,x}(A; (B\cup \{\alpha\})\setminus\{b_1\})),\]
\[\M_{X,Y,z}(G-\{u,s\})=\M_{X,Y,z+2}(S'_{N-1,x+1}(T(A); T(B\setminus\{b_1\}))),\]
\[\M_{X,Y,z}(G-\{v,w\})=\M_{X,Y,z}(S'_{N,x}(A\setminus\{a_{\ell}\}; B \cup\{\alpha\})).\]
Denote by $\mathcal{R}_{u,v}, \mathcal{R}_{w,s}, \mathcal{R}_{u,s}, \mathcal{R}_{v,w},$ the four semihexagons on the right-hand sides above.
By Kuo's condensation Lemma \ref{lemkuo}, we have
\begin{align}\label{proofeq2}
\M_{X,Y,z+2}&(\mathcal{R}_{u,w}) \M_{X,Y,z}(\mathcal{R})=\M_{X,Y,z+2}(\mathcal{R}_{u,v})\M_{X,Y,z}(\mathcal{R}_{w,s})+\M_{X,Y,z+2}(\mathcal{R}_{u,s})\M_{X,Y,z}(\mathcal{R}_{v,w}).
\end{align}

Next, we need to verify that the recurrence (\ref{proofeq2}) the recurrence supports the induction by showing that the $p$-parameters of fives regions different from  $\mathcal{R}=S_{N,x}(A, B)$ are strictly less than $N+t$, and that all other five regions admit a tiling if  $\mathcal{R}$ does. Let us compare how the their $p$-parameters with $p=m+n+t$. We denote the $m$-,$n$-,$t$-, $p$-parameters of the region under consideration by $m'$, $n'$, $t'$, and $p'$.

By investigating each of the five regions, we have the following table
\begin{center}
\setlength{\tabcolsep}{12pt}
\begin{tabular}{lcccc}
\hline
\textbf{Region} & $m'$ & $n'$ & $t'$ & $p'-p$ \\ \hline
$\mathcal{R}_{u,w}$ & $m-1$ & $n$ & $t+b_1-\alpha$ & $-1+b_1-\alpha$ \\[4pt]
$\mathcal{R}_{u,v}$ & $m-1$ & $n$ & $t$ & $-1$ \\[4pt]
$\mathcal{R}_{w,s}$ & $m$ & $n$ & $t+b_1-\alpha$ & $b_1-\alpha$ \\[4pt]
$\mathcal{R}_{u,s}$ & $m$ & $n-1$ & $t+b_1-m-1$ & $b_1-m-2$ \\[4pt]
$\mathcal{R}_{v,w}$ & $m-1$ & $n+1$ & $t+m-\alpha$ & $m-\alpha$ \\[4pt]
\hline
\end{tabular}
\end{center}
One readily sees that all number in the last column are negative as $b_1\leq m<\alpha$ when $t>0$.
\bigskip

Next, we verify that if $\mathcal{R}$ has a tiling, then so is $\mathcal{R}_{u,w}$. Then we can divide both sides of (\ref{proofeq2})  by $\M_{X,Y,z}(\mathcal{R}_{u,w})$ to get
 \begin{align}\label{proofeq2z}
\M_{X,Y,z}(\mathcal{R})=\frac{\M_{X,Y,z+2}(\mathcal{R}_{u,v})\M_{X,Y,z}(\mathcal{R}_{w,s})+\M_{X,Y,z+2}(\mathcal{R}_{u,s})\M_{X,Y,z}(\mathcal{R}_{v,w})}{\M_{X,Y,z+2}(\mathcal{R}_{u,w}) }.
\end{align}
This allows the induction hypothesis to apply two five regions on the right-hand side.
Now assume that $(A,B)$ is a compatible pair, it means
\[|A\cap [j]|+|B\cap [j]| \leq j \tag{*}\]
for every integer $j$ with $1\le j\le N=m+n$. We first see that the translation both sets $A$ and $B$ by the same distance preserves the condition (*). Similarly, the removal of elements and the replacement an element by a larger one in $A$ and $B$ also preserve the condition (*). This implies that $(T(A\setminus \{a_{\ell}\}),\  T((B\cup\{\alpha\})\setminus\{b_1\}))$ is also a compatible pair, and the region $\mathcal{R}_{u,w}$ is tileable as desired.


\bigskip

We will show that the right-hand side of (\ref{proofeq1}) satisfies the same recurrence, then  (\ref{proofeq1}) follows by induction principle.
It remains to prove 
\begin{align}\label{proofeq3}
\mathbf{E}\cdot \M_{X,Y,z+2}&(S'_{N-1,y+1}(T(A\setminus \{a_{\ell}\}), \  T((B\cup\{\alpha\})\setminus\{b_1\}))) \M_{X,Y,z}(S'_{N,y}(A, B))=\notag\\
&\mathbf{H}\cdot \M_{X,Y,z+2}(S'_{N-1,y+1}(T(A\setminus \{a_{\ell}\}),\  T(B)))\M_{X,Y,z}(S'_{N,y}(A,\  (B\cup \{\alpha\})\setminus\{b_1\}))\notag\\
&+\mathbf{K}\cdot \M_{X,Y,z+2}(S'_{N-1,y+1}(T(A),\  T(B\setminus\{b_1\})))\M_{X,Y,z}(S'_{N,y}(A\setminus\{a_{\ell}\},\  B \cup\{\alpha\})),
\end{align}
where
\[\mathbf{E}=g_{x+1,y+1,z+2}(T(A\setminus \{a_{\ell}\}),\ T((B\cup\{\alpha\})\setminus\{b_1\}))g_{x,y,z}(A, B),\]
\[\mathbf{H}=g_{x+1,y+1,z+2}(T(A\setminus \{a_{\ell}\}),\ T(B))g_{x,y,z}(A; (B\cup \{\alpha\})\setminus\{b_1\})),\]
\[\mathbf{K}=g_{x+1,y+1,z+2}(T(A),\  T(B\setminus\{b_1\}))g_{x,y,z}(A\setminus\{a_{\ell}\},\  B \cup\{\alpha\})).\]

This can be done by showing that $\mathbf{E}=\mathbf{H}=\mathbf{K}$. Indeed, if $\mathbf{E}=\mathbf{H}=\mathbf{K}$, then (\ref{proofeq3}) is equivalent to
\begin{align}\label{proofeq4}
\M_{X,Y,z+2}&(S'_{N-1,y+1}(T(A\setminus \{a_{\ell}\}),\ T((B\cup\{\alpha\})\setminus\{b_1\}))) \M_{X,Y,z}(S'_{N,y}(A, B))=\notag\\
&\M_{X,Y,z+2}(S'_{N-1,y+1}(T(A\setminus \{a_{\ell}\}),\  T(B)))\M_{X,Y,z}(S'_{N,y}(A; (B\cup \{\alpha\})\setminus\{b_1\}))\notag\\
&+\M_{X,Y,z+2}(S'_{N-1,y+1}(T(A),\  T(B\setminus\{b_1\})))\M_{X,Y,z}(S'_{N,y}(A\setminus\{a_{\ell}\},\  B \cup\{\alpha\})).
\end{align}
Equation ({proofeq4}) follows in the same way as (\ref{proofeq2}), by Lemma (\ref{lemkuo}).

\bigskip

To complete the proof, we will show that $\mathbf{E}=\mathbf{H}=\mathbf{K}$.

First we consider the ratio
\begin{align}
\frac{\mathbf{E}}{\mathbf{H}}=\frac{g_{x+1,y+1,z+2}(T(A\setminus \{a_{\ell}\}),\  T((B\cup\{\alpha\})\setminus\{b_1\}))}{g_{x,y,z}(A,\  (B\cup \{\alpha\})\setminus\{b_1\})} \cdot \frac{ g_{x,y,z}(A, B)}{g_{x+1,y+1,z+2}(T(A\setminus \{a_{\ell}\}),\  T(B))}.
\end{align}

We use formula (\ref{alternativeg}) and identities (\ref{reduceeq1}), (\ref{reduceeq2}), and (\ref{reduceeq3}) to simplify 
\begin{align}
&\frac{g_{x+1,y+1,z+2}(T(A\setminus \{a_{\ell}\}),\ T((B\cup\{\alpha\})\setminus\{b_1\}))}{g_{x,y,z}(A,\  (B\cup \{\alpha\})\setminus\{b_1\}))}\\
&=\frac{F_x(N)F_{y+1}(N-1)}{F_{x+1}(N-1)F_y(N)} \cdot \frac{\prod_{a\in T(A\setminus \{a_{\ell}\})} Q_{x+1,y+1}(a)}{\prod_{a\in A} Q_{x,y}(a)} \cdot \frac{\prod_{b\in T((B\cup\{\alpha\})\setminus\{b_1\})}Q_{x+1,y+1}(b) P_{x+1,y+1,z+2}(b)}{\prod_{b\in (B\cup\{\alpha\})\setminus\{b_1\}}Q_{x,y}(b)P_{x,y,z}(b)} \notag\\
&= Q_{x,y}(1) \cdot  \frac{\prod_{a\in A\setminus\{a_\ell\}} Q_{x+1,y+1}(a-1)}{\prod_{a\in A} Q_{x,y}(a)} \cdot \frac{\prod_{b \in (B\cup\{\alpha\})\setminus\{b_1\}}Q_{x+1,y+1}(b-1) P_{x+1,y+1,z+2}(b-1)}{\prod_{b\in (B\cup\{\alpha\})\setminus\{b_1\}}Q_{x,y}(b)P_{x,y,z}(b)}\notag\\
&=Q_{x,y}(1)\cdot  \frac{1}{Q_{x,y}(a_{\ell})}\cdot 1=\frac{Q_{x,y}(1)}{Q_{x,y}(a_{\ell})}.
\end{align}
Similarly, we have
\begin{align}
&\frac{ g_{x,y,z}(A, B)}{g_{x+1,y+1,z+2}(T(A\setminus \{a_{\ell}\}),\  T(B)))}\\
&=\frac{F_{x+1}(N-1)F_y(N)}{F_x(N)F_{y+1}(N-1)} \cdot \frac{\prod_{a\in A} Q_{x,y}(a)}{\prod_{a\in T(A\setminus \{a_{\ell}\})} Q_{x+1,y+1}(a)} \cdot \frac{\prod_{b\in B}Q_{x,y}(b)P_{x,y,z}(b)} {\prod_{b\in T(B)}Q_{x+1,y+1}(b) P_{x+1,y+1,z+2}(b)}\notag\\
&= \frac{1}{Q_{x,y}(1)}\frac{\prod_{a\in A} Q_{x,y}(a)}{\prod_{a\in A\setminus \{a_{\ell}\}} Q_{x+1,y+1}(a-1)} \cdot \frac{\prod_{b\in B}Q_{x,y}(b)P_{x,y,z}(b)} {\prod_{b\in B}Q_{x+1,y+1}(b-1) P_{x+1,y+1,z+2}(b-1)}\notag\\
&=\frac{1}{Q_{x,y}(1)}\cdot Q_{x,y}(a_{\ell}) \cdot 1=\frac{Q_{x,y}(a_{\ell})}{Q_{x,y}(1)}.
\end{align}
Multiplying the above two identities, we obtain $\frac{\mathbf{E}}{\mathbf{H}}=1$ or $\mathbf{E}=\mathbf{H}$.
Similarly, we will simplify the ratio
\begin{align}
\frac{\mathbf{E}}{\mathbf{K}}=\frac{g_{x+1,y+1,z+2}(T(A\setminus \{a_{\ell}\}),\  T((B\cup\{\alpha\})\setminus\{b_1\}))}{g_{x,y,z}(A\setminus\{a_{\ell}\},\  B \cup\{\alpha\}))}\cdot \frac{ g_{x,y,z}(A, B)}{g_{x+1,y+1,z+2}(T(A),\ T(B\setminus\{b_1\}))}.
\end{align}
Again, by formula (\ref{alternativeg}) and identities (\ref{reduceeq1}), (\ref{reduceeq2}), and (\ref{reduceeq3})
\begin{align}
&\frac{g_{x+1,y+1,z+2}(T(A\setminus \{a_{\ell}\}),\  T((B\cup\{\alpha\})\setminus\{b_1\}))}{g_{x,y,z}(A\setminus\{a_{\ell}\},\ B \cup\{\alpha\}))}\\
&=\frac{F_x(N)F_{y+1}(N-1)}{F_{x+1}(N-1)F_y(N)} \cdot \frac{\prod_{a\in T(A\setminus \{a_{\ell}\})} Q_{x+1,y+1}(a)}{\prod_{a\in A\setminus\{a_{\ell}\}} Q_{x,y}(a)} \cdot \frac{\prod_{b\in T((B\cup\{\alpha\})\setminus\{b_1\})}Q_{x+1,y+1}(b) P_{x+1,y+1,z+2}(b)}{\prod_{b\in B\cup\{\alpha\}}Q_{x,y}(b)P_{x,y,z}(b)} \notag\\
&=Q_{x,y}(1) \cdot \frac{\prod_{a\in A\setminus \{a_{\ell}\}} Q_{x+1,y+1}(a-1)}{\prod_{a\in A\setminus\{a_{\ell}\}} Q_{x,y}(a)} \cdot \frac{\prod_{b\in (B\cup\{\alpha\})\setminus\{b_1\}}Q_{x+1,y+1}(b-1) P_{x+1,y+1,z+2}(b-1)}{\prod_{b\in B\cup\{\alpha\}}Q_{x,y}(b)P_{x,y,z}(b)} \notag\\
&=Q_{x,y}(1)\cdot 1\cdot \frac{1}{Q_{x,y}(b_1)P_{x,y,z}(b_1)}=\frac{Q_{x,y}(1)}{Q_{x,y}(b_1)P_{x,y,z}(b_1)}.
\end{align}
and
\begin{align}
&\frac{ g_{x,y,z}(A, B)}{g_{x+1,y+1,z+2}(T(A),\  T(B\setminus\{b_1\}))}\\
&=\frac{F_{x+1}(N-1)F_y(N)}{F_x(N)F_{y+1}(N-1)} \cdot \frac{\prod_{a\in A} Q_{x,y}(a)}{\prod_{a\in T(A)} Q_{x+1,y+1}(a)} \cdot \frac{\prod_{b\in B}Q_{x,y}(b)P_{x,y,z}(b)} {\prod_{b\in T((B\setminus\{b_1\})}Q_{x+1,y+1}(b) P_{x+1,y+1,z+2}(b)}\notag\\
&=\frac{1}{Q_{x,y}(1)}\cdot  \frac{\prod_{a\in A} Q_{x,y}(a)}{\prod_{a\in A} Q_{x+1,y+1}(a-1)} \cdot \frac{\prod_{b\in B}Q_{x,y}(b)P_{x,y,z}(b)} {\prod_{b\in B\setminus\{b_1\}}Q_{x+1,y+1}(b-1) P_{x+1,y+1,z+2}(b-1)}\notag\\
&=\frac{1}{Q_{x,y}(1)} \cdot 1 \cdot Q_{x,y}(b_1)P_{x,y,z}(b_1)=\frac{Q_{x,y}(b_1)P_{x,y,z}(b_1)}{Q_{x,y}(1)}.
\end{align}
Multiplying the above two identities we have $\frac{\mathbf{E}}{\mathbf{K}}=1$, or $\mathbf{E}=\mathbf{K}$.
This implies that $\mathbf{E}=\mathbf{H}=\mathbf{K}$ as desired.
\end{proof}

\bibliographystyle{plain}

\bibliography{TriRevision.bib}

\end{document}